\documentclass{amsart}
\usepackage{amssymb}
\usepackage{eucal}
\theoremstyle{plain}
\newtheorem{theorem}{Theorem}

\newtheorem{lemma}{Lemma}

\theoremstyle{definition}

\theoremstyle{remark}

\newtheorem{remark}{Remark}
\numberwithin{equation}{section}
\newcommand{\e}{\epsilon}

\newcommand{\R}{\mathbb R}

\newcommand{\Z}{\mathbb Z}
\newcommand{\C}{\mathbb C}

\newcommand{\Rn}{\mathbb R^n}

\newcommand{\Rm}{\mathbb R^{n+1}}
\begin{document}

\title[Hardy's Uncertainty Principle, Convexity\dots]{Hardy's Uncertainty Principle, Convexity and Schr\"odinger Evolutions}
\author{L. Escauriaza}
\address[L. Escauriaza]{UPV/EHU\\Dpto. de Matem\'aticas\\Apto. 644, 48080 Bilbao, Spain.}
\email{luis.escauriaza@ehu.es}
\thanks{The first and fourth authors are supported  by MEC grant, MTM2004-03029, the second and third authors by NSF grants DMS-0456583 and DMS-0456833 respectively}
\author{C. E. Kenig}
\address[C. E. Kenig]{Department of Mathematics\\University of Chicago\\Chicago, Il. 60637 \\USA.}
\email{cek@math.uchicago.edu}
\author{G. Ponce}
\address[G. Ponce]{Department of Mathematics\\
University of California\\
Santa Barbara, CA 93106\\
USA.}
\email{ponce@math.ucsb.edu}
\author{L. Vega}
\address[L. Vega]{UPV/EHU\\Dpto. de Matem\'aticas\\Apto. 644, 48080 Bilbao, Spain.}
\email{mtpvegol@lg.ehu.es}
\begin{abstract}
We prove the logarithmic convexity of certain quantities, which measure the  quadratic exponential decay at infinity and within two characteristic hyperplanes of solutions of Schr\"odinger evolutions. As a consequence we obtain some uniqueness results that generalize (a weak form of) Hardy's version of the uncertainty principle. We also obtain corresponding results for heat evolutions.
\end{abstract}
\maketitle
\begin{section}{Introduction}\label{S: Introduction} 
In this paper we continue the study initiated in \cite{kpv02} and \cite{ekpv06} on unique continuation properties of solutions of Schr\"odinger evolutions
\begin{equation}
\label{E: 1.1}
i \partial_t u + \triangle u = V(x,t) u\ ,\ \text{in}\  \R^n\times [0,1].
\end{equation}

The goal is to obtain sufficient conditions on a solution $u$, the potential $V$ and the  behavior of the solution at two different times, $t_0=0$ and $t_1=1$, which guarantee that $u\equiv 0$ in $\Rn\times [0,1]$.

One of our motivations  comes from a well known result due to G. H. Hardy \cite[pp. 131]{StSh} (see also \cite{bonamie} for a recent survey on this topic), which concerns the decay of a function $f$ and its Fourier transform, 
\[\hat f(\xi)=(2\pi)^{-\frac n2}\int_{\Rn}e^{-i\xi\cdot x}f(x)\,dx,\]

\emph{If 
$f(x)=O(e^{-|x|^2/\beta^2})$, $\hat f(\xi)=O(e^{-4|\xi|^2/\alpha^2})$ and  $\alpha\beta<4$, then $f\equiv 0$. Also, if $\alpha\beta=4$, $f$ is a constant multiple of $e^{-|x|^2/\beta^2}$,}

\vspace{0,2 cm}

This result can be rewritten in terms of the free solution of the Schr\"odinger equation in $\Rn\times (0,+\infty)$, $i\partial_tu+\triangle u=0$, with initial data $f$,
\begin{equation*}
u(x,t)= (4\pi it)^{-\frac n 2} \int_{\Rn}e^{\frac{ i |x-y|^2}{4t}}f(y)\, dy= \left(2\pi it\right)^{-\frac n2}e^{\frac{i|x|^2}{4t}}\widehat{e^{\frac{i|\,\cdot\,|^2}{4t}}f}\left(\frac x{2t}\right)
\end{equation*}
in the following way: 
\vspace{0,2 cm}

\emph{If $u(x,0)=O(e^{-|x|^2/\beta^2})$, $u(x,T)=O(e^{-|x|^2/\alpha^2})$ and  $\alpha\beta<4T$, then $u\equiv 0$. Also, if $\alpha\beta=4T$, $u$ has as initial data a constant multiple of $e^{-\left(1/\beta^2+i/4T\right)|y|^2}$.}
\vspace{0,2 cm}

The corresponding result in terms of $L^2$-norms and established in \cite{SiSu} is the following: 

\vspace{0,2 cm}
\emph{If $e^{|x|^2/\beta^2}f$, $e^{4|\xi |^2/\alpha^2}\widehat f$ are in $L^2(\Rn)$ and $\alpha\beta\le 4$, then $f\equiv 0$.}

\vspace{0,2 cm}
\emph{If $e^{|x|^2/\beta^2}u(x,0)$, $e^{|\xi |^2/\alpha^2}u(x,T)$ are in $L^2(\Rn)$ and $\alpha\beta\le 4T$, then $u\equiv 0$.}
\vspace{0,2 cm}

In our previous paper \cite {ekpv06} we proved a uniqueness result in this direction for potentials which satisfy
\begin{equation}\label{E: condicionV}
\lim_{R\rightarrow +\infty}\int_0^1\|V(t)\|_{L^\infty(\Rn\setminus B_R)}\,dt =0.
\end{equation}
More precisely, we prove that the only solution to \eqref{E: 1.1} in $C([0,1], H^2(\Rn))$, which verifies that it and its gradient decay faster than any quadratic exponential at times $0$ and $1$ is the zero solution, when $V$ is bounded in $\Rn\times [0,1]$, \eqref{E: condicionV} holds and $\nabla_xV$  is in $L^1_tL^\infty_x(\Rn\times [0,1])$. This linear result was then applied to show that  two regular solutions $u_1$ and $u_2$ of non-linear equations of the type
\begin{equation}
\label{E: NLS}
i\partial_tu+\triangle u=F(u,\overline u), \ \text{in}\ \Rn\times [0,1]
\end{equation}
and for very general non-linearities $F$, must agree in $\Rn\times [0,1]$, when $u_1-u_2$ and its gradient  decay faster than any quadratic exponential at times $0$ and $1$. This replaced the assumption that the solutions coincide on large sub-domains of $\Rn$ at two different times, which was previously studied in \cite{kpv02, Ioke04} and showed that weaker variants of Hardy's Theorem hold even in the context of non-linear Schr\"odinger evolutions.

Our main result in this paper is the following one.

\begin{theorem}\label{T: hardytimeindepent}
Assume that $u$ in $C([0,1]),L^2(\Rn))$ verifies
\begin{equation*}
\label{E: 1.111}\partial_tu=i\left(\triangle u+V(x,t)u\right),\ \text{in}\ \Rn\times [0,1],
\end{equation*}
$\alpha$ and $\beta$ are positive, $\alpha\beta < 2$, $\|e^{\frac{|x|^2}{\beta^2}}u(0)\|_{L^2(\Rn)}$ and $\|e^{\frac{|x|^2}{\alpha^2}}u(1)\|_{L^2(\Rn)}$ are both finite, the potential $V$ is  bounded and either, $V(x,t)=V_1(x)+V_2(x,t)$, with $V_1$ real-valued and \[\sup_{[0,1]}\|e^{\frac{|x|^2}{\left(\alpha t+\beta\left(1-t\right)\right)^2}}V_2(t)\|_{L^\infty(\Rn)}<+\infty\]
or $\lim_{R\rightarrow +\infty}\|V\|_{L^1([0,1], L^\infty(\Rn\setminus B_R)}=0$. Then, $u\equiv 0$.
\end{theorem}

As a direct  consequence of  Theorem \ref{T: hardytimeindepent} we get the following straightforward  result concerning the uniqueness of solutions for non-linear equations of the form \eqref{E: NLS}.

\begin{theorem}
\label{Theorem 1...}

Let $u_1$ and $u_2$ be $C([0,1],H^k(\R^n))$ strong solutions of the equation (1.3)  with $k\in \Z^+$, $k>n/2$,
$F:\C^2\to \C$, $F\in C^k$  and $F(0)=\partial_uF(0)=\partial_{\bar u}F(0)=0$. If there are $\alpha$ and $\beta$ positive  with $\alpha \beta<2$ such that
\begin{equation*}
e^{\frac{|x|^2}{\beta^2}}\left(u_1(0)-u_2(0)\right)\ \text{and}\ e^{\frac{|x|^2}{\alpha^2}}\left(u_1(1)-u_2(1)\right) 
\end{equation*}
are in $L^2(\Rn)$, then $u_1\equiv u_2$.

\end{theorem}

Notice that the condition, $\alpha\beta<2$, is independent of the size of the potential or the dimension and that we do not assume any decay of the gradient neither of the solutions or of time-independent potentials or any regularity of the potentials.  

Our improvement for the results of \cite{ekpv06} comes from  a better understanding of the solutions to \eqref{E: 1.1}, which have a Gaussian decay. We started the study of this particular type of solutions in our recent work \cite{ekpv07}, where we consider free waves  (i.e. $V(x,t)=0$ in \eqref{E: 1.1}) and  among other results  we proved the following:

Assume  that $u$ in $C([0,1], L^2(\Rn))$ is a solution of 
\[\partial_tu-i\triangle u=0,\ \text{in}\ \Rn\times [0,1],\]
and that $\|e^{\gamma |x|^2}u(0)\|$, $\|e^{\gamma |x|^2}u(1)\|$ are both finite. Set $f=e^{\gamma |x|^2}u$ and $H(t)=\left(f,f\right)$. Then, $\log{H(t)}$ is a convex function.
 
The proof of Theorem \ref{T: hardytimeindepent} relies first on extending the above convexity properties to the non-free case, and secondly on a modification of  the definition of the function $H$  as follows: for $e_1=(1,0,\dots , 0)$ and $R>0$ set
\begin{equation}\label{E: 1.1111}
f=e^{\mu|x+Re_1t(1-t)|^2}u,
\end{equation} 
when $0<\mu<\gamma$ and $H(t)=\left(f,f\right)$. Then it is easy to prove at a formal level that

 \[\partial^2_t\log{H(t)}\ge -\frac{R^2}{4\mu}.\]
Therefore
$H(t)e^{-\frac{R^2t(1-t)}{8\mu}}$ is logarithmically convex in $[0,1]$ and
\begin{equation*}
H(t)\le H(0)^{1-t}H(1)^t e^{\tfrac{R^2t(1-t)}{8\mu}}.
\end{equation*}
Taking $t=\frac 12$ and letting $\mu$ increase towards $\gamma$, we have
\begin{equation*}
\int e^{2\gamma |x+\tfrac {Re_1}4|^2}|u(\tfrac 12)|^2\,dx\le \|e^{\gamma |x|^2}u(0)\| \|e^{\gamma |x|^2}u(1)\|e^{\tfrac{R^2}{32\gamma}}.
\end{equation*}
Thus,
\begin{equation*}
\int_{B_{\tfrac{\epsilon R}4}} |u(\tfrac 12)|^2\,dx\le \|e^{\gamma |x|^2}u(0)\| \|e^{\gamma |x|^2}u(1)\|e^{\tfrac{R^2\left(1-4\gamma^2(1-\epsilon)^2\right)}{32\gamma}},
\end{equation*}
when $0<\e<1$, which implies that $u\equiv 0$ by letting $R$ tend to infinity, when $\gamma > \frac 12$.

The path that goes from the formal level to a rigorous one is not an easy one. In fact in section \ref{S:  A positive commutator and a misleading frequency function} we will give explicit examples of functions $H(t)$ such that $\log H$ is formally convex and however the corresponding inequalities lead  to false statements. Therefore most of this paper is devoted to make rigorous the above argument. The starting point is to prove similar properties to those obtained in \cite{ekpv08} for free solutions. One of the results we get is the following one.

\begin{theorem}\label{T: densityfortimedependentpotentialsvega}
Assume that $u$ in $C([0,1]),L^2(\Rn))$ verifies 
\begin{equation*}
\partial_tu=i\left(\triangle u+V(x,t)u\right),\ \text{in}\ \Rn\times [0,1],
\end{equation*}
$V=V_1(x)+V_2(x,t)$,  $V_1$ is real-valued, $\|V_1\|_{\infty}\le M_1$ and that there are positive numbers $\alpha$ and $\beta$ such that
\[\|e^{\frac{|x|^2}{\beta^2}}u(0)\|\ ,\ \|e^{\frac{|x|^2}{\alpha^2}}u(1)\|\ \text{and}\ \sup_{[0,1]}\|e^{\frac{|x|^2}{\left(\alpha t+\left(1-t\right)\beta\right)^2}}V_2(t)\|_{\infty}<+\infty.\]
Then, $\|e^{\frac{|x|^2}{\left(\alpha t+\left(1-t\right)\beta\right)^2}}u(t)\|^{\alpha t+\left(1-t\right)\beta}$ is \lq\lq logarithmically convex\rq\rq\  in $[0,1]$ and there is $N=N(\alpha,\beta)$ such that
\begin{equation*}
\|e^{\frac{|x|^2}{\left(\alpha t+\left(1-t\right)\beta\right)^2}}u(t)\|\\\le e^{N\left(M_1+M_2+M_1^2+M_2^2\right)}\|e^{\frac{|x|^2}{\beta^2}} u(0)\|^{\frac{\beta\left(1-t\right)}{\alpha t+\beta\left(1-t\right)}}\|e^{\frac{|x|^2}{\alpha^2}} u(1)\|^{\frac{\alpha t}{\alpha t+\beta\left(1-t\right)}},
\end{equation*}
when $0\le s\le 1$ and $M_2=\sup_{[0,1]}{\| e^{\frac{|x|^2}{\left(\alpha t+\beta\left(1-t\right)\right)^2}} V_2(t)\|_{\infty}}\, e^{2\sup_{[0,1]}\|\Im V_2(t)\|_{\infty}}$. Moreover,
\begin{equation*}
\begin{split}
\|\sqrt{t(1-t)}\, &e^{\frac{|x|^2}{\left(\alpha t+\left(1-t\right)\beta\right)^2}}\nabla  u\|_{L^2(\Rn\times [0,1])}\\&\le Ne^{N\left(M_1+M_2+M_1^2+M_2^2\right)}
\left[\|e^{\frac{|x|^2}{\beta^2}}u(0)\|+\|e^{\frac{|x|^2}{\alpha^2}}u(1)\|\right].
\end{split}
\end{equation*}
\end{theorem}

In order to prove this theorem we have to approximate the solution using some artificial diffusion. The corresponding results are interesting in themselves and can be found in section \ref{S:  Algunos Lemmas}. As a byproduct we get  examples of solutions to \eqref{E: 1.1} which have Gaussian decay, when the potential $V$ is time independent. It is enough to consider as initial data the solution, at  say time one, of the corresponding heat equation that at time zero is a Gaussian. This property was already established in \cite{ekpv08} for free solutions, and it turned out to be a characterization of those Gaussian solutions. It would be interesting to prove similar characterizations for variable coefficient Hamiltonians. Also in section \ref{S:  Algunos Lemmas} we give an abstract result, Lemma \ref {L: freq1}, that shows how to get logarithmic convexity properties from the positivity of some specific commutators. It turns out that these commutators are the same as the ones that appear in the proof of the $L^2$-Carleman estimates we used in our previous paper \cite{ekpv06}.
In fact, the weight $\mu|x+Re_1t(1-t)|^2$ that appears in \eqref{E: 1.1111} is a refinement of the ones used in \cite{ekpv06}.

We are indebted to E. Zuazua for pointing out the following application of Hardy's uncertainty principle to prove the following optimal decay result for solutions of the free heat equation (See also \cite[Section 5]{j06}):

	{\it If $f$  and $e^{\frac{|x|^2}{\delta^2}}e^\triangle f$ are in $L^2(\Rn)$ for some $\delta\le 2$. Then, $f\equiv 0$.}
	
	In fact, applying Hardy's uncertainty principle to $e^\triangle f$, $e^{\frac{|x|^2}{\delta^2}}e^\triangle f$ and $e^{\frac{4|\xi|^2}{2^2}}\widehat{e^\triangle f}=\widehat f$ would be in $L^2(\Rn)$, and $2\delta\le 4$ implies $e^{\triangle}f\equiv 0$. Then,  backward uniqueness arguments, see for instance  \cite[Chapter 3, Theorem 11]{lm60} or \cite[Chapter 3]{e98}, show that $f\equiv 0$. Here, we prove the following weaker extension of this result for parabolic operators with variable coefficientes.
	
\begin{theorem}\label{T: toremaparabolicco}
Let $u$ in $L^\infty([0,1]), L^2(\Rn))\cap L^2([0,T], H^1(\Rn))$ verify
\begin{equation*}
\begin{cases}
\partial_tu=\triangle u+V(x,t)u,\ \text{in}\ \Rn\times (0,1],\\
u(0)=f,
\end{cases}
\end{equation*}
where $V$ is bounded in $\Rn\times [0,1]$ and assume that $f$ and $e^{\frac{|x|^2}{\delta^2}}u(1)$ are in $L^2(\Rn)$ for some $\delta <1$. Then, $f\equiv 0$ in $\Rn$.
\end{theorem}

		It is natural to expect  that the Hardy uncertainty principle holds  on Schr\"odinger and heat evolutions with bounded potentials and with parameters $\alpha$, $\beta$ or $\delta$ verifing the condition of the free case.
		
		In the sketch of the proof of Theorem \ref{T: hardytimeindepent}  that we have done above we have assumed that $\alpha=\beta$. That one can easily reduce to this case is proved in section \ref {S: 3} using the so called conformal transformation or Appell transform. In section \ref{S:  Variable Coefficients. } we prove Theorem \ref{T: densityfortimedependentpotentialsvega}, in section \ref{S:  A Hardy Type Uncertainty Principle}  we give the proof of Theorem \ref{T: hardytimeindepent}, in section \ref{S:  A positive commutator and a misleading frequency function} we give some examples of some misleading convex functions and in section \ref{S: Parabolic analog} we prove Theorem \ref{T: toremaparabolicco}.

\end{section}
\begin{section}{A few Lemmas}\label{S:  Algunos Lemmas}
In the sequel \[\left(f,g\right)=\int_{\Rn}f\overline g\,dx\ ,\ \|f\|^2=\left(f,f\right)\ ,\ f^+=\max{\{f,0\}}\]
and $\|f\|_{\infty}$ denotes the $L^\infty$-norm of $f$ over $\Rn$.
\begin{lemma}\label{L: parabolicdecay}
Assume that  $u$ in $L^\infty([0,1], L^2(\Rn))\cap L^2([0,1]),H^1(\Rn))$ satisfies
\begin{equation*}
\partial_tu=(A+iB)\left(\triangle u+V(x,t)u+F(x,t)\right),\ \text{in}\ \Rn\times (0,1],\\
\end{equation*}
$A>0$ and $B\in \R$. Then,
\begin{multline*}
e^{-M_T}\|e^{\frac{\gamma A|x|^2}{A+4\gamma\left(A^2+B^2\right)T}}u(T)\|\\\le \|e^{\gamma |x|^2}u(0)\|+\sqrt{A^2+B^2}\|e^{\frac{\gamma A|x|^2}{A+4\gamma\left(A^2+B^2\right)t}}F(t)\|_{L^1([0,T], L^2(\Rn))},
\end{multline*}
when $\gamma\ge 0$, $0\le T\le 1$ and $M_T=\|A\left(\text{\it Re\,}V\right)^+-B\text{\it Im\,}V\|_{L^1([0,T],L^\infty(\Rn))}$.
\end{lemma}
\begin{proof}
Write $v=e^\varphi u$, where $\varphi$ is a real-valued function to be chosen later. The function $v$ verifies
\[\partial_tv=\mathcal Sv+\mathcal Av+\left(A+iB\right)e^{\varphi}F,\  \text{in}\  \Rn\times (0,1]\ ,\]
where the symmetric and skew-symmetric operators $\mathcal S$ and $\mathcal A$ are given by
\begin{align*}
\mathcal S=&A\left(\triangle +|\nabla\varphi|^2\right)-iB\left(2\nabla\varphi\cdot\nabla+\triangle\varphi\right)+\left(\partial_t\varphi+A\text{\it Re\,}V-B\text{\it Im\,}V\right)\ ,\\
\mathcal A= &iB\left(\triangle +|\nabla\varphi|^2\right)-A\left(2\nabla\varphi\cdot\nabla+\triangle\varphi\right)+ i\left(B\text{\it Re\,}V+A\text{\it Im\,}V\right)\ .
\end{align*}

To prove Lemma \ref{L: parabolicdecay} we use the energy method and try to keep track of the decay of the $L^2(\Rn)$-norm of $v$. Formally,
\begin{equation*}
\partial_t\|v\|^2=2\text{\it Re\,}\left(\mathcal Sv,v\right)+2\text{\it Re\,}\left(\left(A+iB\right)e^{\varphi}F,v\right)\ ,
\end{equation*}
when $t\ge 0$. Again, a formal integration by parts gives that
\begin{multline}\label{E: f—rmula m‡gica}
\text{\it Re\,}(\mathcal Sv,v)=-A\int_{\Rn}|\nabla v|^2\,dx+\int_{\Rn}\left(A|\nabla\varphi|^2+\partial_t\varphi\right)|v|^2\,dx\\+2B\,\text{Im}\int_{\Rn}\overline v\nabla\varphi\cdot\nabla v\,dx+\int_{\Rn}\left(A\text{\it Re\,}V-B\text{\it Im\,}V\right)|v|^2\,dx
\end{multline}
and the Cauchy-Schwarz's inequality implies that
\begin{align*}
\partial_t\|v(t)\|^2\le & \,2 \|A\left(\text{\it Re\,}V(t)\right)^+-B\text{\it Im\,}V(t)\|_{\infty}\|v(t)\|^2\\
+ &2\sqrt{A^2+B^2}\|e^{\varphi} F(t)\|\|v(t)\|\ ,
\end{align*}
when 
\begin{equation}\label{E: desigualdad necesaria}
\left(A+\frac{B^2}{A}\right)|\nabla\varphi|^2+\partial_t\varphi\le 0,\  \text{in}\ \Rm_+ .
\end{equation}

When $\varphi(x,t)= a(t)\phi(x)$, it suffices that
\begin{equation*}
a^2(t)\left(A+\frac{B^2}{A}\right)|\nabla \phi(x)|^2+a'(t)\phi(x)\le 0.
\end{equation*}
At the end we shall require  that $\phi(x)= |x|^2$. In that case the latter holds, when
\begin{equation}\label{E: ecuaci—n de a}
\begin{cases}
a'(t)=-4\left(A+\frac{B^2}{A}\right)a^2(t),\\
a(0)=\gamma .
\end{cases}
\end{equation}

To formalize the integration by parts and calculations carried out above, given $\gamma >0$, we truncate $|x|^2$ as
\begin{equation*}
\phi_R(x)=\begin{cases}
 |x|^2,\ &|x|\le R,\\
 R^2,\ &|x|>R,
 \end{cases}
 \end{equation*}
 regularize $\phi_R$ with a radial mollifier $\theta_\rho$ and set 
 \[\varphi_{\rho, R}(x,t)=a(t)\,\theta_\rho\ast\phi_{R}(x)\ ,\ v_{\rho, R}=e^{\varphi_{\rho, R}}u,\] 
 where
\begin{equation*}
a(t)=\frac{\gamma A}{A+4\gamma\left(A^2+B^2\right)t}\ .
\end{equation*} 
is  the solution to \eqref{E: ecuaci—n de a}.

Because the right hand side of \eqref{E: f—rmula m‡gica} only involves the first derivatives of $\varphi$, $\phi_R$ is Lipschitz and bounded at infinity, 
\begin{equation*}
\theta_\rho\ast\phi_R\le \theta_\rho\ast |x|^2=|x|^2+C(n)\rho^2 
\end{equation*}
and \eqref{E: desigualdad necesaria} holds uniformly in the variables $\rho$ and $R$, when $\varphi$ is replaced by $\varphi_{\rho,R}$, it follows (and now rigorously) that the estimate
\begin{equation*}
\|v_{\rho,R}(T)\|\le e^{M_T}\left(\|e^{\gamma |x|^2}u(0)\|+\sqrt{A^2+B^2}\|e^{\varphi_{\rho, R}}F\|_{L^1([0,T], L^2(\Rn))}\right)
\end{equation*}
holds uniformly in $\rho$ and $R$. Lemma \ref{L: parabolicdecay} follows after letting $\rho$ tend to zero and $R$ to  infinity. 
\end{proof}
\begin{lemma}\label{L: freq1}
$\mathcal S$ is a symmetric operator, $\mathcal A$ is skew-symmetric, both are allowed to depend on the time variable, $G$ is a positive function, $f(x,t)$ is a reasonable function, 
\begin{equation*}
H(t)=\left( f, f\right)\ ,\ D(t)=\left( \mathcal Sf, f\right)\ ,\ \partial_t\mathcal S=\mathcal S_t
\quad \text{and}\quad N(t)=\frac{D(t)}{H(t)}\ .
\end{equation*}
Then,
\begin{multline}
\label{E: derivadasegunda}
\partial_t^2H= 2\partial_t\text{\it Re}\left(\partial_tf-\mathcal Sf-\mathcal Af,f\right)+2\left(\mathcal S_tf+\left[\mathcal S,\mathcal A\right]f,f\right)\\+\|\partial_tf-\mathcal Af+\mathcal Sf\|^2-\|\partial_tf-\mathcal Af-\mathcal Sf\|^2
\end{multline}
and
\begin{equation*}
\dot N(t)\ge \left(\mathcal S_tf +\left[\mathcal S,\mathcal A\right]f, f\right)/H- \|\partial_tf-\mathcal Af-\mathcal Sf\|^2/\left(2H\right).
\end{equation*}
Moreover, if
\begin{equation}\label{E: condicionesbase}
|\partial_tf-\mathcal Af-\mathcal Sf|\le M_1|f| +G,\ \text{in}\ \Rn\times [0,1],\quad  \mathcal S_t+\left[\mathcal S,\mathcal A\right]\ge -M_0,
\end{equation}
and 
\[M_2=\sup_{[0,1]}{\|G(t)\|/\|f(t)\|}\]
is finite, then 
$\log H(t)$ is  \lq\lq logarithmically convex\rq\rq\ in $[0,1]$ and there is a universal constant $N$ such that
\begin{equation}\label{E: convexidadlogaritmica}
H(t)\le e^{N\left(M_0+M_1+M_2+M_1^2+M_2^2\right)}H(0)^{1-t}H(1)^t,\ \text{when}\ 0\le t\le 1.
\end{equation}
\end{lemma}
\begin{proof} Formally,
\begin{align*}
\dot H(t)&=2\text{\it Re}\left(\partial_tf,f\right)=2\text{\it Re}\left(\partial_tf-\mathcal Sf-\mathcal Af,f\right)+2\left(\mathcal Sf,f\right) 
\end{align*}
and
\begin{equation}\label{E: logaritmo1}
\dot H(t)=2\text{\it Re}\left(\partial_tf-\mathcal Sf-\mathcal Af,f\right)+2D(t)\ .
\end{equation}
Also,
\begin{align*}
\dot H(t)&=\text{\it Re}\left(\partial_tf+\mathcal Sf,f\right)+\text{\it Re}\left(\partial_tf-\mathcal Sf,f\right),\\
D(t)&=\frac 12\text{\it Re}\left(\partial_tf+\mathcal Sf,f\right)-\frac 12\text{\it Re}\left(\partial_tf-\mathcal Sf,f\right)\
\end{align*}
and multiplying the last two formulae,
\[
\dot H(t)D(t)=\frac 12\left(\text{\it Re}\left(\partial_tf+\mathcal Sf,f\right)\right)^2-\frac 12\left(\text{\it Re}\left(\partial_tf-\mathcal Sf,f\right)\right)^2.  
\]

Adding an antisymmetric operator does not change the real parts, and so
\begin{equation}\label{E: logaritmo3}
\dot H(t)D(t)=\frac 12\left(\text{\it Re}\left(\partial_tf-\mathcal Af+\mathcal Sf,f\right)\right)^2-\frac 12\left(\text{\it Re}\left(\partial_tf-\mathcal Af-\mathcal Sf,f\right)\right)^2\ .  
\end{equation}

Differentiating $D(t)$,
\begin{align*}
\dot D(t)&=\left(\mathcal S_tf,f\right)+\left(\mathcal S\partial_tf,f\right)+\left(\mathcal Sf,\partial_tf\right)\\
&=\left(\mathcal S_tf,f\right)+2\text{\it Re}\left(\partial_tf,\mathcal Sf\right)\\
&=\left(\mathcal S_tf+\left[\mathcal S,\mathcal A\right]f,f\right)+2\text{\it Re}\left(\partial_tf-\mathcal Af,\mathcal Sf\right)  
\end{align*}
and the polarization identity gives 
\begin{equation}\label{E: logaritmo2}
\dot D(t)= \left(\mathcal S_tf+\left[\mathcal S,\mathcal A\right]f,f\right)+\frac 12\|\partial_tf-\mathcal Af+\mathcal Sf\|^2-\frac 12\|\partial_tf-\mathcal Af-\mathcal Sf\|^2.
\end{equation}
The formula  \eqref{E: derivadasegunda} for the second derivative of $H$ follows from \eqref{E: logaritmo1} and \eqref{E: logaritmo2}. 
The identity
\begin{align*}
\dot N(t)&= \left(\mathcal S_tf +\left[\mathcal S,\mathcal A\right]f, f\right)/H\\
& +\frac 12\left[\|\partial_tf-\mathcal Af+\mathcal Sf\|^2\|f\|^2-\left(\text{\it Re}\left(\partial_tf-\mathcal Af+\mathcal Sf,f\right)\right)^2\right]/H^2\\
&+ \frac 12\left[ \left(\text{\it Re}\left(\partial_tf-\mathcal Af-\mathcal Sf,f\right)\right)^2-\|\partial_tf-\mathcal Af-\mathcal Sf\|^2\|f\|^2\right]/H^2
\end{align*}
follows from \eqref{E: logaritmo3} and \eqref{E: logaritmo2}. The inequality  in Lemma \ref{L: freq1} follows from the positiveness of the second line (Cauchy-Schwarz's inequality) and of the fourth term on the right hand side of the previous identity.

When \eqref{E: condicionesbase} holds, the first part in Lemma \ref{L: freq1} shows  that
\[\dot N(t)\ge -\left( M_0+M_1^2+M_2^2\right),\]
and from \eqref{E: logaritmo1}
\[\partial_t\left[\log{H(t)}+\text{O}\left(1\right)\right]=2N(t).\]

All together,
\[\partial^2_t\left(\log{H(t)}+\text{O}\left(1\right)\right)\ge 0,\  \text{when}\ 0\le t\le 1,\]
where $\text{O}\left(1\right)$ is a function verifying, $|\text{O}\left(1\right)|\le N\left(M_0+M_1+M_2+M_1^2+M_2^2\right)$ in $[0,1]$. The integration of the inequality
\begin{equation*}
\partial_s\left(\log{H(s)}+\text{O}\left(1\right)\right)\le\partial_\tau\left(\log{H(\tau)}+\text{O}\left(1\right)\right),\ \text{when}\ 0\le s\le t\le\tau\le 1,
\end{equation*}
over the intervals, $0\le s\le t$ and $t\le\tau\le 1$, implies \eqref{E: convexidadlogaritmica}.
\end{proof}

\begin{lemma}\label{L: convexlogsch}
Assume that $u$ in $L^\infty([0,1]),L^2(\Rn))\cap L^2([0,1],H^1(\Rn))$ verifies 
\begin{equation}\label{E: initial value sch}
\partial_tu=(A+iB)\left(\triangle u+V(x,t)u+F(x,t)\right),\ \text{in}\ \Rn\times [0,1],
\end{equation}
where $A>0$, $B\in\R$, $V$ is complex-valued, $\gamma >0$ and $\sup_{[0,1]}\|V(t)\|_{\infty}\le M_1$. Set 
\[M_2=\sup_{[0,1]}{\|e^{\gamma |x|^2}F(t)\|/\|u(t)\|}\]
and assume that  $\|e^{\gamma |x|^2}u(0)\|$, $\|e^{\gamma |x|^2}u(1)\|$ and $M_2$ are finite. Then, 
$\|e^{\gamma|x|^2}u(t)\|$ is \lq\lq logarithmically convex\rq\rq\  in $[0,1]$ and there is a universal constant $N$ such that
\begin{multline}\label{E: resultadologaritmico}
\|e^{\gamma|x|^2}u(t)\|\le\\ e^{N\left[\left(A^2+B^2\right)\left(\gamma M_1^2+M_2^2\right)+\sqrt{A^2+B^2}\left(M_1+M_2\right)\right]}\|e^{\gamma|x|^2}u(0)\|^{1-t}\|e^{\gamma|x|^2}u(1)\|^t,
\end{multline}
when $0\le t\le 1$.
\end{lemma}

\begin{proof}
Let $f=e^{\gamma\varphi} u$, where $\varphi=\varphi(x,t)$ is to be chosen. The function $f$ verifies
\begin{equation}\label{E: loquefcumple}
\partial_tf=\mathcal Sf+\mathcal Af +(A+iB)\left(Vf+e^{\gamma\varphi}F\right),\ \text{in}\ \Rm_+,
\end{equation}
with symmetric and skew-symmetric operators $\mathcal S$ and $\mathcal A$
\begin{equation}
\label{E: formulaoperadores}
\begin{aligned}
\mathcal S=&A\left(\triangle +\gamma^2|\nabla\varphi|^2\right)-iB\gamma\left(2\nabla\varphi\cdot\nabla+\triangle\varphi\right)+\gamma\partial_t\varphi,\\
\mathcal A= &iB\left(\triangle +\gamma^2|\nabla\varphi|^2\right)-A\gamma\left(2\nabla\varphi\cdot\nabla+\triangle\varphi\right).
\end{aligned}
\end{equation}

A calculation shows that,
\begin{multline}\label{E: formulaconmutadorindependientetiempo}
\mathcal S_t+\left[\mathcal S,\mathcal A\right]=\gamma\partial^2_t\varphi+4\gamma^2A\nabla\varphi\cdot\nabla\partial_t\varphi-2iB\gamma\left(2\nabla\partial_t\varphi\cdot\nabla+\triangle\partial_t\varphi\right)\\ -\gamma\left(A^2+B^2\right)\left[4\nabla\cdot\left(D^2\varphi\nabla\,\,\,\right)-4\gamma^2D^2\varphi\nabla\varphi\cdot\nabla\varphi+\triangle^2\varphi\right].
\end{multline}

At the end we shall require that $\varphi(x,t)= |x|^2$, where
\begin{equation*}
\mathcal S_t+\left[\mathcal S,\mathcal A\right]= -\gamma\left(A^2+B^2\right)\left[8\triangle -32\gamma^2|x|^2\right]
\end{equation*}
and
\begin{equation}\label{E: conmutadormasderivada}
\left(\mathcal S_tf+\left[\mathcal S,\mathcal A\right]f,f\right)=\gamma\left(A^2+B^2\right)
\int_{\Rn}8|\nabla f|^2+32\gamma^2|x|^2|f|^2\,dx\ .
\end{equation}
This identity, the boundedness of $V$ and \eqref{E: loquefcumple} imply that
\begin{equation}
\label{E: loqueciumplenfyconmutador}
|\partial_tf-\mathcal Sf-\mathcal Af|\le \sqrt{A^2+B^2}\left(M_1|f|+e^{\gamma\varphi}|F|\right)\quad ,\quad \mathcal S_t+\left[\mathcal S,\mathcal A\right]\ge 0,
\end{equation}
and if we knew that the quantities and calculations involved in the proof of Lemma \ref{L: freq1}  were  finite and correct, when $f=e^{\gamma |x|^2}u$, we would have the \lq\lq logarithmic  convexity\rq\rq\  of $H(t)=\|e^{\gamma |x|^2}u(t)\|^2$ and  get \eqref{E: resultadologaritmico} from Lemma \ref{L: freq1}.

To justify the validity of the previous arguments, given $a$ and $\rho$ in $(0,1)$, define
\begin{equation*}
\varphi_a(x) =
\begin{cases}
|x|^2,\  & |x|<1\ ,\\
\left(2|x|^{2-a}-a\right)/(2-a),\ & |x|\ge 1
\end{cases}
\end{equation*}
and replace $\varphi= |x|^2$ by $\varphi_{a,\rho}=\theta_\rho\ast\varphi_a$, where $\theta$ in $C_0^\infty(\Rn)$ is a radial function. Observe that $\varphi_a$ is a $C^{1,1}(\Rn)$ convex function, $\varphi_{a,\rho}\le |x|^2+C(n)\rho^2$, $\varphi_{a,\rho}$ is convex and grows at infinity not faster than $|x|^{2-a}$. At the same time,
\begin{equation*}
\triangle\varphi_a(x)=
\begin{cases}
2n,\ & |x|\le 1,\\
2(n-a)|x|^{-a},\ &|x|\ge 1,
\end{cases}
\end{equation*}
and the distribution $\partial_j\triangle\varphi_a$, $j=1,\dots,n$, is equal to
\begin{equation*}
-2ax_j\, d\sigma -2a(n-a)x_j|x|^{-a-2}\chi_{\Rn\setminus B_1},
\end{equation*}
where $d\sigma$ is surface measure on $\partial B_1$. This and the identity 
\begin{equation*}
\triangle^2\varphi_{a,\rho}=\sum_{j=1}^n\partial_j\theta_\rho\ast\partial_j\triangle\varphi_a,
\end{equation*}
show that 
\begin{equation}\label{E: tama–odelbilaplaciano}
\|\triangle^2\varphi_{a,\rho}\|_{\infty}\le C(n,\rho)a. 
\end{equation} 
Set then, $f_{a,\rho}=e^{\gamma\varphi_{a,\rho}}u$ and $H_{a,\rho}(t)=\|f_{a,\rho}\|^2$ in Lemma \ref{L: freq1}. The decay bound in Lemma \ref{L: parabolicdecay} and the interior regularity for solutions of \eqref{E: initial value sch} (Here we use that $A$ is positive) can now be used qualitatively to make sure that the quantities or calculations involved in the proof of Lemma \ref{L: freq1} are finite and correct for $f_{a,\rho}$. In this case, $f_{a,\rho}$ verifies
\begin{equation}\label{E: loquefcumple2}
\partial_tf_{a,\rho}=\mathcal S^{a,\rho}f_{a,\rho}+\mathcal A^{a,\rho}f_{a,\rho} +(A+iB)\left(Vf_{a,\rho}+e^{\gamma\varphi_{a,\rho}}F\right),\ \text{in}\ \Rn\times [0,1],
\end{equation}
with symmetric and skew-symmetric operators $\mathcal S^{a,\rho}$ and $\mathcal A^{a,\rho}$ given by \eqref{E: formulaoperadores} with  $\varphi$ replaced by $\varphi_{a,\rho}$. The formula for the operator
\[\mathcal S^{a,\rho}_t+\left[\mathcal S^{a,\rho},\mathcal A^{a,\rho}\right]\]
 in \eqref{E: formulaconmutadorindependientetiempo}, the convexity of $\varphi_{a,\rho}$, the bounds \eqref{E: tama–odelbilaplaciano} and \eqref{E: loquefcumple2} imply that the inequalities
\begin{align*}
|\partial_tf_{a,\rho}-\mathcal S^{a,\rho}f_{a,\rho}-\mathcal A^{a,\rho}f_{a,\rho}|&\le \sqrt{A^2+B^2}\left(M_1|f_{a,\rho}|+e^{\gamma\varphi_{a,\rho}}F\right),\\ S^{a,\rho}_t+\left[\mathcal S^{a,\rho}, \mathcal A^{a,\rho}\right] &\ge 0,
\end{align*}
hold and $M_2(a,\rho)\le e^{C(n)\rho^2} M_2$, when $0< a,\, \rho\le 1$. In particular, $H_{a, \rho}$ is \lq\lq logarithmically convex\rq\rq\ in $[0,1]$ and
\begin{equation}\label{E: convexidadlogarritmica3}
H_{a,\rho}(t)\le e^{N\left[\left(A^2+B^2\right)\left(M_1^2+M_2^2\right)+\sqrt{A^2+B^2}\left(M_1+M_2\right)\right]}H_{a,\rho}(0)^{1-t}H_{a,\rho}(1)^t.
\end{equation}

Then, \eqref{E: resultadologaritmico} follows  after taking first  the limit, when $a$ tends to zero in \eqref{E: convexidadlogarritmica3} and then, when $\rho$ tends to zero.
\end{proof}

\begin{lemma}\label{L: regularidad}
Assume that $A+iB$, $u$ and $V$ are as  in Lemma \ref{L: convexlogsch} and $\gamma>0$. Then, 
\begin{multline}\label{E: controlgradiente}
\|\sqrt{t(1-t)}e^{\gamma|x|^2}\nabla u\|_{L^2(\Rn\times [0,1])}+\|\sqrt{t(1-t)}|x|e^{\gamma|x|^2}u\|_{L^2(\Rn\times [0,1])}\\ \le N\left[\left(1+M_1\right)\sup_{[0,1]} \|e^{\gamma |x|^2}u(t)\|+\sup_{[0,1]}\| e^{\gamma |x|^2}F\|_{L^2(\Rn\times [0,1])}\right],
\end{multline}
where $N$ remains bounded, when $\gamma$ and $A^2+B^2$ are bounded  below.
\end{lemma}

\begin{proof}
A formal integration by parts shows that
\begin{equation*}
\int_{\Rn}|\nabla f|^2+4\gamma^2|x|^2|f|^2\,dx=\int_{\Rn}e^{2\gamma |x|^2}\left(|\nabla u|^2-2n\gamma |u|^2\right)\,dx,
\end{equation*}
when $f=e^{\gamma |x|^2}u$, while either well known properties of Hermite functions \cite{Than} or integration by parts, the Cauchy-Schwarz's inequality and the identity, $n=\nabla\cdot x$, give that 
\begin{equation*}
\int_{\Rn}|\nabla f|^2+4\gamma^2|x|^2|f|^2\,dx\ge 2\gamma n\int_{\Rn}|f|^2\,dx.
\end{equation*}
The sum of the last two formulae gives the inequality
\begin{equation}\label{E: desigualdadinteresante}
2\int_{\Rn}|\nabla f|^2+4\gamma^2|x|^2|f|^2\,dx\ge\int_{\Rn}e^{2\gamma |x|^2}|\nabla u|^2\,dx.
\end{equation}

Integration over $[0,1]$ of $t(1-t)$ times the formula \eqref{E: derivadasegunda} for the second derivative of $H(t)=\|f(t)\|^2$ and integration by parts, shows that in the general framework of Lemma \ref{L: freq1}
\begin{multline}\label{E: laderivadasegundadeH}
2\int_0^1t(1-t)\left(\mathcal S_tf+\left[\mathcal S,\mathcal A\right]f,f\right)dt+2\int_0^1H(t)\,dt\le H(1)+H(0)\\+2\int_0^1(1-2t)\text{\it Re}\left(\partial_tf-\mathcal Sf-\mathcal Af,f\right)\,dt + \int_0^1t(1-t)\|\partial_tf-\mathcal Af-\mathcal Sf\|^2\,dt.
\end{multline}

Assuming again that the last two calculations are justified for $f=e^{\gamma |x|^2}u$, \eqref{E: laderivadasegundadeH}, \eqref{E: conmutadormasderivada}, \eqref{E: loqueciumplenfyconmutador}, \eqref{E: desigualdadinteresante} and the identity, $\nabla f=e^{\gamma |x|^2}\left(\nabla u+2\gamma x u\right)$, imply the Lemma.

The interior regularity of  the solutions to \eqref{E: initial value sch} (Here we use again that $A>0$) shows that the calculations leading to \eqref{E: desigualdadinteresante} and \eqref{E: laderivadasegundadeH} are justified, when $f=e^{\left(\gamma-\rho\right)|x|^2}u$, $0<\rho<\gamma$,  and the right hand side of \eqref{E: controlgradiente} is finite. The Lemma follows letting $\rho$ tend to zero.
\end{proof}
\section{The conformal or Appell transformation}\label{S: 3}
\begin{lemma}\label{L: transformada} Assume that  $u(y,s)$ verifies
\begin{equation*}
\partial_su=(A+Bi)\left(\triangle u+V(y,s)u+F(y,s)\right)\ ,\ \text{in}\  \R^n\times [0,1],
\end{equation*}
$A+iB\neq 0$, $\alpha$ and $\beta$ are positive, $\gamma\in\R$ and set
\[\widetilde u(x,t)=\left(\tfrac{\sqrt{\alpha\beta}}{\alpha(1-t)+\beta t}\right)^{\frac n2}u\left(\tfrac{\sqrt{\alpha\beta}\, x}{\alpha(1-t)+\beta t}, \tfrac{\beta t}{\alpha(1-t)+\beta t}\right)e^{\frac{\left(\alpha-\beta\right) |x|^2}{4(A+iB)(\alpha(1-t)+\beta t)}}.\]
Then, $\widetilde u$ verifies
 \begin{equation*}
\partial_t\widetilde u=(A+Bi)\left(\triangle \widetilde u+\widetilde V(x,t)\widetilde u+\widetilde F(x,t)\right)\ ,\ \text{in}\  \R^n\times [0,1],
\end{equation*}
with
\begin{equation*}
\begin{split}
&\widetilde V(x,t)=\tfrac{\alpha\beta}{\left(\alpha(1-t)+\beta t\right)^2}\,V\left(\tfrac{\sqrt{\alpha\beta}\, x}{\alpha(1-t)+\beta t}, \tfrac{\beta t}{\alpha(1-t)+\beta t}\right),\\
\widetilde F(x,t)&=\left(\tfrac{\sqrt{\alpha\beta}}{\alpha(1-t)+\beta t}\right)^{\frac n2+2}F\left(\tfrac{\sqrt{\alpha\beta}\, x}{\alpha(1-t)+\beta t}, \tfrac{\beta t}{\alpha(1-t)+\beta t}\right)e^{\frac{\left(\alpha-\beta\right) |x|^2}{4(A+iB)(\alpha(1-t)+\beta t)}}.
\end{split}
\end{equation*}
Moreover,
\begin{equation*}
\|e^{\gamma |x|^2}\widetilde F(t)\|=\tfrac{\alpha\beta}{\left(\alpha\left(1-t\right)+\beta t\right)^2}\|e^{\left[\frac{\gamma\alpha\beta}{(\alpha s+\beta(1-s))^2}+\frac{\left(\alpha-\beta\right) A}{4(A^2+B^2)(\alpha s+\beta(1-s))}\right]|y|^2}F(s)\|
\end{equation*}
and
\begin{equation*}
\|e^{\gamma |x|^2}\widetilde u(t)\| = \|e^{\left[\frac{\gamma\alpha\beta}{(\alpha s+\beta(1-s))^2}+\frac{\left(\alpha-\beta\right) A}{4(A^2+B^2)(\alpha s+\beta(1-s))}\right]|y|^2}u(s)\|,
\end{equation*}
when $s=\tfrac{\beta t}{\alpha(1-t)+\beta t}$ and $\gamma\in\R$.
\end{lemma}

\begin{proof} 
When $u$ satisfies
\begin{equation}\label{E: parab—licageneral}
\partial_su=(A+Bi)\left(\triangle u+H(y,s)\right)\ ,\ \text{in}\  \R^n\times [0,1],
\end{equation} 
the function, $u_1(x,t)=u(\sqrt rx,rt+\tau)$, verifies 
\[\partial_t u_1= \left(A+iB\right)\left(\triangle u_1+ r H(\sqrt rx , rt+\tau)\right)\]
and $u_2(x,t)=t^{-\frac n2}u(x/t,1/t)e^{\frac{|x|^2}{4\left(A+iB\right)t}}$ is a solution to 
\[\partial_tu_2=-\left(A+Bi\right)\left(\triangle u_2+t^{-\frac n2-2}H(x/t,1/t)e^{\frac{|x|^2}{4\left(A+iB\right)t}}\right).\]   
These two facts and the sequel of changes of variables below prove the Lemma, when $\alpha>\beta$ :
\[u\left(\sqrt{\tfrac{\alpha\beta}{\alpha-\beta}}\, x , \tfrac{\alpha\beta}{\alpha-\beta}\, t-\tfrac\beta{\alpha-\beta}\right)\]
is a solution to the same non-homogeneous equation but with right-hand side
\[\tfrac{\alpha\beta}{\alpha-\beta} H\left(\sqrt{\tfrac{\alpha\beta}{\alpha-\beta}}\,x,\tfrac{\alpha\beta}{\alpha-\beta}\, t-\tfrac \beta{\alpha-\beta}\right).\] 
The function,
\[\tfrac{1}{(\alpha-t)^{\frac n2}}\,u\left(\tfrac{\sqrt{\alpha\beta}\, x}{\sqrt{\alpha-\beta} (\alpha-t)}, \tfrac{\alpha\beta}{\left(\alpha-\beta\right)(\alpha-t)}-\tfrac\beta{\alpha-\beta}\right)e^{\frac{|x|^2}{4(A+iB)(\alpha-t)}}\]
verifies \eqref{E: parab—licageneral} with right-hand side 
 \[\tfrac{\alpha\beta}{\left(\alpha-\beta\right) (\alpha -t)^{\frac n2+2}}H\left(\tfrac{\sqrt{\alpha\beta}\, x}{\sqrt{\alpha-\beta} (\alpha-t)}, \tfrac{\alpha\beta}{\left(\alpha-\beta\right)(\alpha-t)}-\tfrac\beta{\alpha-\beta}\right)e^{\frac{|x|^2}{4(A+iB)(\alpha-t)}}.\]
Replacing $(x,t)$ by $(\sqrt{\alpha-\beta}\, x,\left(\alpha-\beta\right) t)$,
\begin{equation}\label{E: casiterminado}
\tfrac {1}{(\alpha(1-t)+\beta t)^{\frac n2}}\,u\left(\tfrac{\sqrt{\alpha\beta}\, x}{\alpha(1-t)+\beta t}, \tfrac{\alpha\beta}{\left(\alpha-\beta\right)(\alpha(1-t)+\beta t)}-\tfrac\beta{\alpha-\beta}\right)e^{\frac{\left(\alpha-\beta\right) |x|^2}{4(A+iB)(\alpha(1-t)+\beta t)}}
\end{equation}
is a solution to \eqref{E: parab—licageneral} but with right-hand
\begin{equation}\label{E: aladerecha}
\tfrac{\alpha\beta}{(\alpha+\beta -\alpha t)^{\frac n2+2}}H\left(\tfrac{\sqrt{\alpha\beta}\, x}{\alpha(1-t)+\beta t}, \tfrac{\alpha\beta}{\left(\alpha-\beta\right)(\alpha(1-t)+\beta t)}-\tfrac\beta{\alpha-\beta}\right)e^{\frac{\left(\alpha-\beta\right) |x|^2}{4(A+iB)(\alpha(1-t)+\beta t)}}.
\end{equation}
Finally, observe that
\[s= \tfrac{\beta t}{\alpha(1-t)+\beta t}= \tfrac{\alpha\beta}{\left(\alpha-\beta\right)(\alpha(1-t)+\beta t)}-\tfrac\beta{\alpha-\beta}\]
and multiply \eqref{E: casiterminado} and \eqref{E: aladerecha}  by $\left(\sqrt{\alpha\beta}\right)^{\frac n2}$. 

The case $\beta>\alpha$ follows by reversing the time with the changes of variables, $s'=1-s$ and $t'=1-t$. The relations between the different norms of $\widetilde u$, $u$, $\widetilde F$ and $F$ follow undoing the changes of variables and using the identity
\begin{equation*}
\tfrac{\sqrt{\alpha\beta}}{\alpha\left(1-t\right)+\beta t}=\tfrac{\alpha s+\beta\left(1-s\right)}{\sqrt{\alpha\beta}}\,.
\end{equation*}
\end{proof}
\end{section}
\begin{section}{Variable Coefficients. Proof of Theorem \ref{T: densityfortimedependentpotentialsvega}.}\label{S:  Variable Coefficients. }
We are ready to prove Theorem \ref{T: densityfortimedependentpotentialsvega}.
\begin{proof}  We may assume that $\alpha \neq\beta$. The case $\alpha=\beta$ follows from the latter by replacing $\beta$ by $\beta+\delta$, $\delta>0$, and letting $\delta$ tend to zero. We may also assume that $\alpha<\beta$. Otherwise, replace $u$ by $\overline u(1-t)$. Set then, $H=\triangle+V_1(x)$ and let $e^{t\left(A+iB\right)H}u_0$ denote the $C([0,1],L^2(\Rn))$-solution to
\begin{equation*}
\begin{cases}
\partial_tv=\left(A+iB\right)\left(\triangle v+V_1(x)v\right),\ \text{in}\ \Rn\times [0,1],\\
v(0)=u_0,
\end{cases}
\end{equation*}
when $\Re \left(A+iB\right) \ge 0$. The Duhamel principle shows that
\begin{equation}\label{E: representacion}
u(t)=e^{itH}u(0)+i\int_0^t e^{i\left(t-s\right)H}\left(V_2(s)u(s)\right)\,ds\ ,\ \text{in}\ \Rn\times [0,1].
\end{equation}
For $0\le\e \le 1$, set
\begin{equation}\label{E: formulamagica1}
F_\e(t)=\tfrac i{\e+i}e^{\e tH}\left(V_2(t)u(t)\right),
\end{equation}
and
\begin{equation}\label{E: formulamagica2}
u_\e (t)=e^{\left(\e +i\right)tH}u(0)+\left(\e+i\right)\int_0^t e^{\left(\e+i\right)\left(t-s\right)H}F_\e(s)\,ds.
\end{equation}
Then, $u_\e$ is in $L^\infty([0,1],L^2(\Rn))\cap L^2([0,1], H^1(\Rn))$ and verifies
\begin{equation*}
\begin{cases}
\partial_t u_\e=\left(\e+i\right)\left(H u_\e+F_\e(t)\right),\ \text{in}\ \Rn\times [0,1],\\
u_\e(0)=u(0).
\end{cases}
\end{equation*}
The identities \cite{pazy}
 \begin{equation}\label{E: semigrupos}
e^{\left(z_1+z_2\right)H}=e^{\left(z_2+z_1\right)H}=e^{z_1H}e^{z_2H}\ ,\ \text{when}\  \Re z_1, \Re z_2\ge 0,
\end{equation}
\eqref{E: representacion}, \eqref{E: formulamagica1} and \eqref{E: formulamagica2} show that
\begin{equation}\label{E: unaformulamassenilla}
u_\e(t)=e^{\e tH}u(t)\ , \text{when} \ 0\le t\le 1.
\end{equation}
In particular,
\[u_\e(1)=e^{\epsilon H}u(1)\] 
and Lemma \ref{L: parabolicdecay} with $A+iB=\e$, $\gamma =\frac 1{\beta^2}$, $F\equiv 0$ and the fact that $u_\e(0)=u(0)$ imply that
\begin{equation*}
\|e^{\frac {|x|^2}{\beta^2+4\e}}u_\e (1)\|\le e^{\e \|V_1\|_{L^\infty(\Rn)}}\|e^{\frac{|x|^2}{\beta^2}}u(1)\|\quad ,\quad \|e^{\frac{|x|^2}{\alpha^2}}u_\epsilon(0)\|= \|e^{\frac{|x|^2}{\alpha^2}}u(0)\|.
\end{equation*}
A second application of Lemma \ref{L: parabolicdecay} with $A+iB=\epsilon$, $F\equiv 0$, the value of $\gamma=\frac 1{\left(\alpha t+\beta\left(1-t\right)\right)^2}$ and \eqref{E: formulamagica1} show that
\begin{equation*}
\|\e^{\frac{|x|^2}{\left(\alpha t+\beta\left(1-t\right)\right)^2+4\e t}}F_\e(t)\|\le e^{\epsilon \|V_1\|_{\infty}}\| \e^{\frac{|x|^2}{\left(\alpha t+\beta\left(1-t\right)\right)^2}} V_2(t)\|_{\infty}\|u(t)\|,
\end{equation*}
when $0\le t\le 1$. Setting, $\alpha_\e=\alpha+2\e$ and $\beta_\e=\beta+2\e$, the last three inequalities give that
\begin{equation}\label{E: loquecumpleuepsilon}
\| e^{\frac {|x|^2}{\beta_e^2}}u_\e (1)\|\le e^{\e \|V_1\|_{\infty}}\|e^{\frac{|x|^2}{\beta^2}}u(1)\|\quad ,\quad \|e^{\frac{|x|^2}{\alpha_\e^2}}u_\epsilon(0)\|\le \|e^{\frac{|x|^2}{\alpha^2}}u(0)\|,
\end{equation}
\begin{equation}\label{E: loquecumpleF}
\|e^{\frac{|x|^2}{\left(\alpha_\e t+\beta_\e\left(1-t\right)\right)^2}}F_\e(t)\|\le e^{\epsilon \|V_1\|_{\infty}}\| e^{\frac{|x|^2}{\left(\alpha t+\beta\left(1-t\right)\right)^2}} V_2(t)\|_{\infty}\|u(t)\|.
\end{equation}
A third application of Lemma \ref{L: parabolicdecay} with  $A+iB=\epsilon$, $F\equiv 0$, $\gamma =0$, and \eqref{E: formulamagica1}, \eqref{E: unaformulamassenilla} implies that
\begin{equation}\label{E: masdesigualdades}
\|F_\e(t)\|\le e^{\epsilon \|V_1\|_{\infty}}\|V_2(t)\|_{L^\infty(\Rn)}\|u(t)\|\ ,\ \|u_\e(t)\|\le e^{\epsilon \|V_1\|_{\infty}}\|u(t)\|,
\end{equation}
when $0\le t\le 1$. Set then, $\gamma_\e=\frac 1{\alpha_\e\beta_\e}\,$ and let
\begin{equation*}
\widetilde u_\epsilon (x,t)= \left(\tfrac{\sqrt{\alpha_\epsilon \beta_\e}}{\alpha_\epsilon (1-t)+\beta_\e t}\right)^{\frac n2}u_\epsilon \left(\tfrac{\sqrt{\alpha_\epsilon \beta_\e}\, x}{\alpha_\epsilon (1-t)+\beta_\e t}, \tfrac{\beta_\e t}{\alpha_\epsilon (1-t)+\beta_\e t}\right)e^{\frac{\left(\alpha_\epsilon-\beta_\e\right) |x|^2}{4(\epsilon+i)(\alpha_\epsilon(1-t)+\beta_\e t)}},
\end{equation*}
be the function associated to $u_\e$ in Lemma \ref{L: transformada}, when $A+iB=\epsilon+i$ and $\alpha$, $\beta$ are replaced respectively by $\alpha_\epsilon$ and $\beta_\e$. Because $\alpha<\beta$, $\widetilde u_\e$ is in $L^\infty([0,1]),L^2(\Rn))\cap L^2([0,1],H^1(\Rn))$ and satisfies
\begin{equation*}
\partial_t \widetilde u_\e=\left(\e+i\right)\left(\triangle \widetilde u_\e+\widetilde V^\e_1(x,t) \widetilde u_\e+\widetilde F_\e(t)\right),\ \text{in}\ \Rn\times [0,1],
\end{equation*}
where $\widetilde V_1^\e$ is real-valued,
\begin{equation}\label{E: otraformulahorrorosa}
\widetilde V_1^\e(x,t)=\tfrac{\alpha_\e \beta_\e}{\left(\alpha_\e (1-t)+\beta_\e t\right)^2}\,V_1\left(\tfrac{\sqrt{\alpha_\e \beta_\e}\, x}{\alpha_\e (1-t)+\beta_\e t}\right)\quad ,\quad \sup_{[0,1]}{\|\widetilde V_1^\e(t)\|_\infty}\le\tfrac{\beta}\alpha M_1,
\end{equation}
\begin{equation*}
\widetilde F_\e(x,t)=\left(\tfrac{\sqrt{\alpha_\e\beta_\e}}{\alpha_\e(1-t)+\beta_\e t}\right)^{\frac n2+2}F_\e\left(\tfrac{\sqrt{\alpha_\e\beta_\e}\, x}{\alpha_\e(1-t)+\beta_\e t}, \tfrac{\beta_\e t}{\alpha_\e(1-t)+\beta_\e t}\right)e^{\frac{\left(\alpha_\e-\beta_\e\right) |x|^2}{4\left(\e+i\right)\left(\alpha_\e(1-t)+\beta_\e t\right)}},
\end{equation*}
\begin{equation}\label{E: estohacelatransformadaconF}
\|e^{\gamma_\e |x|^2}\widetilde F_\e(t) \|\le \tfrac\beta\alpha\,\|\e^{\frac{|x|^2}{\left(\alpha_\e t+\beta_\e\left(1-t\right)\right)^2}}F_\e(s)\|\ ,\  \|\widetilde F_\e(t) \|\le \tfrac\beta\alpha\|F_\e(s)\|,
\end{equation}
and
\begin{equation}\label{E: estohacelatransformadaconue}
\begin{split}
\|e^{\gamma_\e |x|^2}\widetilde u_\e(t)\| = \|e&^{\left[\frac{1}{(\alpha_\e s +\beta_\e (1-s))^2}+\frac{\left(\alpha_\e-\beta_\e\right) A}{4(A^2+B^2)(\alpha_\e s+\beta_\e(1-s))}\right]|y|^2}u_\e(s)\|,\\
 &\|\widetilde u_\e(t)\| \le \|u_\e(s)\|,
\end{split}
\end{equation}
when $s=\tfrac{\beta_\e t}{\alpha_\e(1-t)+\beta_\e t}$. The above identity, when $t$  is  zero or one  and  \eqref{E: loquecumpleuepsilon} shows that
\begin{equation}\label{E: condicionesextremos}
\|e^{\gamma_\e |x|^2}\widetilde u_\e(0)\|\le \|e^{\frac{|x|^2}{\beta^2}}u(0)\|\ ,\ \|e^{\gamma_\e |x|^2}\widetilde u_\e(1)\|\le e^{\e\|V^+\|_{\infty}}\|e^{\frac{|x|^2}{\beta^2}}u(1)\|.
\end{equation}
On the other hand,
\begin{equation}\label{E: algoinmediatosobreu}
N_1^{-1}\|u(0)\|\le\|u(t)\|\le N_1\|u(0)\|\ ,\ \text{when}\ 0\le t\le 1\ ,\ N_1= e^{\sup_{[0,1]}\|\Im V_2(t)\|_{\infty}},
\end{equation}
and the equation satisfied by $\widetilde u_\e$ and the energy method imply that
\begin{equation}\label{E: unadesigualdad}
\partial_t\|\widetilde u_\e(t)\|^2\le 2\e \|\widetilde V^\e_1(t)\|_{\infty}\|\widetilde u_\e(t)\|^2+2\|\widetilde F_\e(t)\|\|\widetilde u_\e(t)\|.
\end{equation}
Let, $0=t_0<t_1<t_2<\dots<t_m=1$, be a uniformly distributed partition of $[0,1]$, where $m$ will be chosen later. The inequality \eqref{E: unadesigualdad}, \eqref{E: otraformulahorrorosa}, the inequality in \eqref{E: estohacelatransformadaconue}, the second inequality in \eqref{E: estohacelatransformadaconF}, \eqref{E: masdesigualdades} and \eqref{E: algoinmediatosobreu} imply that there is $N_2$, which depends on $\frac{\beta}{\alpha}$, $\|V_1\|_{L^{\infty}(\Rn)}$ and $\sup_{[0,1]}\|V_2(t)\|_{\infty}$, such that
\begin{equation}\label{E: energia}
\|\widetilde u_\e(t_i)\|\le e^{\tfrac{\e\beta}\alpha\|V_1\|_{\infty}}\|\widetilde u_\e(t)\|+N_2\sqrt{t_i-t_{i-1}}\|u(0)\|,
\end{equation}
when $t_{i-1}\le t\le t_i$, $0<\e \le 1$ and $i=1,\dots ,m$. Choose now $m$ so that
\begin{equation}\label{E: limita}
N_2\max_{1\le i\le m}\sqrt{t_i-t_{i-1}}\le \frac1{4N_1},
\end{equation}
where $N_1$ was defined in \eqref{E: algoinmediatosobreu}. Because, $\lim_{\e \rightarrow 0^+}\|\widetilde u_\e (t)\|= \|u(s)\|$, when $s=\tfrac{\beta t}{\alpha(1-t)+\beta t}$ and \eqref{E: algoinmediatosobreu}, there is $\e_0$ such that
\begin{equation}\label{E: otradesigualdad}
\|\widetilde u_\e(t_i)\|\ge \frac 1{2N_1}\,\|u(0)\|\ ,\  \text{when}\ \ 0<\e\le \e_0\ ,\ i=1,\dots ,m,
\end{equation}
and now, \eqref{E: otradesigualdad}, \eqref{E: limita} and \eqref{E: energia} show that
\begin{equation}\label{E: loquedebiaocurrir}
\|\widetilde u_\e(t)\| \ge \frac 1{4N_1}\|u(0)\|\ ,\ \text{when}\ 0<\e\le\e_0\ ,\ 0\le t\le 1.
\end{equation}
It is now simple to verify that  \eqref{E: loquedebiaocurrir}, the first inequality in \eqref{E: estohacelatransformadaconF}, \eqref{E: loquecumpleF} and \eqref{E: algoinmediatosobreu} imply that 
\begin{equation}\label{E: loquecuestaescribir}
\sup_{[0,1]}{\frac{\|e^{\gamma_\e |x|^2}\widetilde F_\e(t)\|}{\|\widetilde u_\e (t)\|}}\le
\tfrac{4\beta}\alpha\  M_2(\e)\ ,\ \text{when}\ 0<\e\le\e_0,
\end{equation}
where
\begin{equation*}
M_2(\e)=e^{2\sup_{[0,1]}\|\Im V_2(t)\|_{\infty}+\e\|V_1\|_{\infty}}\sup_{[0,1]}{\| e^{\frac{|x|^2}{\left(\alpha t+\beta\left(1-t\right)\right)^2}} V_2(t)\|_{\infty}}.
\end{equation*}
We can use Lemma \ref{L: convexlogsch}, \eqref{E: condicionesextremos}, \eqref {E: otraformulahorrorosa}
 and \eqref{E: loquecuestaescribir} to show that $\|e^{\gamma_\e |x|^2}\widetilde u_\e (t)\|$ is \lq\lq logarithmically convex\rq\rq\  in $[0,1]$ and that
\begin{equation}\label{E: casillegando}
\|e^{\gamma_\e|x|^2}\widetilde u_\e(t)\|\le e^{N\left(M_1+M_2(\e)+M_1^2+M_2(\e)^2\right)}\|e^{\frac{|x|^2}{\beta^2}}u(0)\|^{1-t}\|e^{\frac{|x|^2}{\alpha^2}}u(1)\|^t,
\end{equation}
when $0\le t\le 1$ and $0<\e\le\e_0$ and with $N=N(\alpha,\beta)$. Then, Lemma \ref{L: regularidad} gives that
\begin{equation*}
\begin{split}
\|\sqrt{t(1-t)}\,e^{\gamma _\e |x|^2}\nabla \widetilde u_\e\|_{L^2(\Rn\times  [0 ,1])}+ & \|\sqrt{t(1-t)}\,|x|\,e^{\gamma _\e |x|^2}\nabla\widetilde u_\e\|_{L^2(\Rn\times [0,1])}\\\le Ne^{N\left(M_1+M_2(\e)+M_1^2+M_2(\e)^2\right)} & \left[\|e^{\frac{|x|^2}{\beta^2}}u(0)\|+\|e^{\frac{|x|^2}{\alpha^2}}u(1)\|\right],
\end{split}
\end{equation*}
when $0<\e\le\e_0$, and the \lq\lq logarithmic convexity\rq\rq\  and  regularity of $u$ follow from the limit of the identity in \eqref{E: estohacelatransformadaconue}, the final limit relation between the variables $s$ and $t$, $s=\tfrac{\beta t}{\alpha(1-t)+\beta t}$, and letting $\e$ tend to zero in \eqref{E: casillegando} and the above inequality. 
\end{proof}
\begin{remark} We thank R. Killip for pointing out the following application of Lemma \ref{L: parabolicdecay} and the identities \eqref{E: semigrupos}  to generate Gaussian decaying solutions of $\partial_t=iH$, when $H=\triangle+V_1(x)$ and $V_1$ verify the conditions in Theorem \ref{T: densityfortimedependentpotentialsvega}. In fact, if $e^{\gamma |x|^2}u_0$ is in $L^2(\Rn)$ and $u(t)=e^{itH}\left(e^Hu_0\right)$, we have $u(t)=e^{\left(\frac 1t +i\right)tH}u_0$ and from Lemma \ref{L: parabolicdecay} 
\[\|e^{\frac {\gamma|x|^2}{1+4\gamma\left(1+t^2\right)}}u(t)\|\le e^{\|V_1\|_{\infty}}\|e^{\gamma |x|^2}u_0\|\ ,\ \text{when}\ t\ge 0.\]
\end{remark}

Next, we recall the  following result established in \cite{kpv02}:

\begin{lemma}\label{L: casoL1} There are $N$ and $\epsilon_0>0$ such that the following holds: 

If $\lambda$ is in $\Rn$, $V$ is a complex-valued potential, $\|V\|_{L^1([0,1], L^{\infty}(\Rn))}\leq \epsilon_0$ and $u\in C([0,1],L^2(\mathbb R^n))$ satisfies 
\begin{equation*}
\partial_tu=i\left(\Delta u +V(x,t)u+F(x,t)\right)\ ,\ \text{in}\ \Rn\times [0,1].
\end{equation*}
Then,
\begin{equation*}
\sup_{[0,1]}\| e^{\lambda\cdot x} u(t)\| \leq N
\left[\|e^{\lambda\cdot x} u(0)\| + \|e^{\lambda\cdot x} u(1)\| +
\|e^{\lambda\cdot x} F(t)\|_{L^1([0,1],L^2(\Rn))}\right].
\end{equation*}

\end{lemma}

\begin{theorem}\label{T: elcasoantiguo}
Assume that $u$ in $C([0,1]),L^2(\Rn))$ verifies 
\begin{equation*}
\partial_tu=i\left(\triangle u+V(x,t)u\right),\ \text{in}\ \Rn\times [0,1],
\end{equation*}
where $V$ is in $L^\infty(\Rn\times [0,1])$, $\lim_{R\rightarrow +\infty}\|V\|_{L^1([0,1], L^\infty(\Rn\setminus B_R))}=0$, $\alpha$ and $\beta$ are positive and $\|e^{\frac{|x|^2}{\beta^2}}u(0)\|$, $\|e^{\frac{|x|^2}{\alpha^2}}u(1)\|$ are finite. Then, there is $N=N(\alpha,\beta)$ such that
\begin{multline*}
\sup_{[0,1]}\|e^{\frac{|x|^2}{\left(\alpha t+\left(1-t\right)\beta\right)^2}}u(t)\|+\|\sqrt{t(1-t)}e^{\frac{|x|^2}{\left(\alpha t+\left(1-t\right)\beta\right)^2}}\nabla u\|_{L^2(\Rn\times [0,1])}\\\ \le
Ne^{N\sup_{[0,1]}\|V(t)\|_\infty}\left[\|e^{\frac{|x|^2}{\beta^2}}u(0)\|+\|e^{\frac{|x|^2}{\alpha^2}}u(1)\|+\sup_{[0,1]}\|u(t)\|\right].
\end{multline*}
\end{theorem}

\begin{proof}
Set $\gamma =1/\left(\alpha\beta\right)$ and let
\begin{equation}\label{E: dequeestahecha}
\widetilde u (x,t)= \left(\tfrac{\sqrt{\alpha \beta}}{\alpha (1-t)+\beta t}\right)^{\frac n2}u \left(\tfrac{\sqrt{\alpha\beta}\, x}{\alpha (1-t)+\beta t}, \tfrac{\beta t}{\alpha (1-t)+\beta t}\right)e^{\frac{\left(\alpha -\beta\right) |x|^2}{4i(\alpha (1-t)+\beta t)}}
\end{equation}
denote the function associated in Lemma \ref{L: transformada}  to $u$, when $A+iB=i$. This function is in $C([0,1], L^2(\Rn))$ and verifies
\begin{equation*}
\partial_t \widetilde u=i\left(\triangle\widetilde u+\widetilde V(x,t) \widetilde u\right),\ \text{in}\ \Rn\times [0,1],\\
\end{equation*}
with
\begin{equation*}\label{E: otraformulahorrorosaparaV}
\widetilde V(x,t)=\tfrac{\alpha\beta}{\left(\alpha (1-t)+\beta t\right)^2}\,V\left(\tfrac{\sqrt{\alpha \beta}\, x}{\alpha (1-t)+\beta t}\, , \tfrac{\beta t}{\alpha (1-t)+\beta t}\right),
\end{equation*}
\begin{equation*}
\sup_{[0,1]}\|\widetilde V(t)\|_\infty\le \max{\{\tfrac\alpha\beta\, ,\tfrac\beta\alpha\}}\sup_{[0,1]}\| V(t)\|_\infty\ ,\ \lim_{R\rightarrow +\infty}\|\widetilde V\|_{L^1([0,1], L^\infty(\Rn\setminus B_R))}=0
\end{equation*}
and	
\begin{equation}
\label{E: estohacelatransformadaconu4}
\begin{aligned}
&\|e^{\gamma |x|^2}\widetilde u(t)\| = \|e^{\frac{|x|^2}{(\alpha s+\beta (1-s))^2}}u(s)\|,\\
&\|\widetilde u(t)\|=\|u(s)\|\ ,\,\; \text{when}\;\;\ s= \frac{\beta t}{\alpha (1-t)+\beta t}\ . 
\end{aligned}
\end{equation}
Choose $R>0$  such that $\|\widetilde V\|_{L^1([0,1],L^\infty(\Rn\setminus B_R))}\le \e_0$. Then,
\begin{equation*}
\partial_t\widetilde u=i\left(\triangle \widetilde u+\widetilde V_R(x,t)\widetilde u+\widetilde F_R(x,t)\right),
\end{equation*}
with $\widetilde V_R(x,t)=\chi_{\Rn\setminus B_R}\widetilde V(x,t)$, $\widetilde F_R=\chi_{B_R}\widetilde V(x,t)\widetilde u$, and using Lemma \ref{L: casoL1}
\begin{equation*}
\sup_{[0,1]}\|e^{\lambda\cdot x}\widetilde u(t)\|\le N\left[\|e^{\lambda\cdot x}\widetilde u(0)\|+\|e^{\lambda\cdot x}\widetilde u(1)\|+e^{|\lambda |R}\sup_{[0,1]}\|\widetilde V(t)\|_{\infty}\sup_{[0,1]}\|\widetilde u(t)\|\right].
\end{equation*} 

Replace $\lambda$ by $\lambda\sqrt\gamma$ in the above inequality, square both sides, multiply all by $e^{-|\lambda|^2/2}$ and integrate both sides with respect to $\lambda$ in $\Rn$. This and the identity,
\begin{equation*}
\int_{\mathbb R^n}e^{2\sqrt\gamma\,\lambda\cdot x-\frac{|\lambda|^2}{2}}\,d\lambda =(2\pi)^{n/2}
e^{2\gamma |x|^2},
\end{equation*}
imply the inequality
\begin{equation*}
\sup_{[0,1]}\|e^{\gamma |x|^2}\widetilde u(t)\|\le N\left[\|e^{\gamma |x|^2}\widetilde u(0)\|+\|e^{\gamma |x|^2}\widetilde u(1)\|+e^{2\gamma R^2}\sup_{[0,1]}\|\widetilde V(t)\|_{\infty}\sup_{[0,1]}\|\widetilde u(t)\|\right].
\end{equation*}
This inequality and \eqref{E: estohacelatransformadaconu4} imply that
\begin{equation}\label{E: unabuenacota}
\sup_{[0,1]}\|e^{\gamma |x|^2}\widetilde u(t)\|\le\\
N\left[\|e^{\frac{|x|^2}{\beta^2}}u(0)\|+\|e^{\frac{|x|^2}{\alpha^2}}u(1)\|+\sup_{[0,1]}\| V(t)\|_{\infty}\sup_{[0,1]}\|u(t)\|\right],
\end{equation}
for some new constant $N$.

To prove the regularity of $u$ we proceed as in \eqref{E: representacion}, \eqref{E: formulamagica1} and \eqref{E: formulamagica2}. The Duhamel formula shows that
\begin{equation}\label{E: representacion3}
\widetilde u(t)=e^{it\triangle}\widetilde u(0)+i\int_0^t e^{i\left(t-s\right)\triangle}\left(\widetilde V(s)\widetilde u(s)\right)\,ds\ ,\ \text{in}\ \Rn\times [0,1].
\end{equation}
For $0<\e< 1$, set
\begin{equation}\label{E: formulamagica13}
\widetilde F_\e(t)=\tfrac i{\e+i}e^{\e t\triangle}\left(\widetilde V(t)\widetilde u(t)\right),
\end{equation}
and
\begin{equation}\label{E: formulamagica23}
\widetilde u_\e (t)=e^{\left(\e +i\right)t\triangle}\widetilde u(0)+\left(\e+i\right)\int_0^t e^{\left(\e+i\right)\left(t-s\right)\triangle }\widetilde F_\e(s)\,ds.
\end{equation}
The identities \cite{pazy}
 \begin{equation*}
e^{\left(z_1+z_2\right)\triangle}=e^{\left(z_2+z_1\right)\triangle}=e^{z_1\triangle}e^{z_2\triangle}\ ,\ \text{when}\  \Re z_1, \Re z_2\ge 0,
\end{equation*}
\eqref{E: representacion3}, \eqref{E: formulamagica13} and \eqref{E: formulamagica23} show that
\begin{equation}\label{E: unaformulamassenilla3}
\widetilde u_\e(t)=e^{\e t\triangle}\widetilde u(t)\ ,\  \text{when}\ 0\le t\le 1,
\end{equation}
and from Lemma \ref{L: parabolicdecay} with $A+iB=\e$, \eqref{E: unaformulamassenilla3} and \eqref{E: formulamagica13},
\begin{equation}\label{E: loqueverificauepsilon}
\begin{split}
&\sup_{[0,1]}\|e^{\gamma_\e |x|^2}\widetilde u_\e (t)\|\le \sup_{[0,1]}\|e^{\gamma |x|^2}\widetilde u(t)\|,\\
\sup_{[0,1]}\|&e^{\gamma_\e |x|^2}\widetilde F_\epsilon(t)\|\le e^{\sup_{[0,1]}\|\widetilde V(t)\|_\infty}\sup_{[0,1]}\|e^{\gamma |x|^2}\widetilde u(t)\|,
\end{split}
\end{equation}
when $\gamma_\e=\frac\gamma{1+4\gamma\e}$. Then, Lemma \ref{L:  regularidad}, \eqref{E: loqueverificauepsilon} and \eqref{E: unabuenacota} show that
\begin{multline*}
\|\sqrt{t(1-t)}\,e^{\gamma_\e |x|^2}\nabla \widetilde u_\e\|_{L^2(\Rn\times [0,1])}+\|\sqrt{t(1-t)}\,|x|\,e^{\gamma_\e |x|^2} \widetilde u_\e\|_{L^2(\Rn\times [0,1])}\\
\le Ne^{N\sup_{[0,1]}\|V(t)\|_\infty} \left[\|e^{\frac{|x|^2}{\beta^2}}u(0)\|+\|e^{\frac{|x|^2}{\alpha^2}}u(1)\|+\sup_{[0,1]}\|u(t)\|\right].
\end{multline*}
The Lemma follows from this inequality, \eqref{E: estohacelatransformadaconu4}, \eqref{E: unabuenacota}, \eqref{E: dequeestahecha} and letting $\e$ tend to zero.
\end{proof}
\end{section}
\begin{section}{A Hardy Type Uncertainty Principle. Proof of Theorem \ref{T: hardytimeindepent}.}\label{S:  A Hardy Type Uncertainty Principle}

As we mentioned in the introduction, the motivation behind the Carleman inequality in Lemma \ref{L: lonecesarioparamejorar} below is the following monotonicity or frequency function argument related to Lemma \ref{L: freq1}: 

When $u$ in $C([0,1], L^2(\Rn))$ is a free solution to the free Schr\"odinger equation
\[\partial_tu-i\triangle u=0,\ \text{in}\ \Rn\times [0,1],\]
$\|e^{\gamma |x|^2}u(0)\|$, $\|e^{\gamma |x|^2}u(1)\|$ are both finite, $f=e^{\mu |x+Rt(1-t)|^2-\frac{R^2t(1-t)}{8\mu}}u$ and $H=\left(f,f\right)$. Then, $\log H$ is logaritmicaly convex in $[0,1]$, when $0<\mu<\gamma$.

The formal application of the above argument to a $C([0,1], L^2(\Rn))$ solution to 
\begin{equation}\label{E: unaecuacion}
\partial_tu-i\left(\triangle u+V(x,t)u\right)=0,\ \text{in}\ \Rn\times [0,1],
\end{equation}
implies a similar result, when $V$ is a bounded potential, though the justification of the correctness of the manipulations involved in the corresponding formal application of Lemma \ref{L: freq1} are not obvious to us. In fact, we can only justify these manipulations, when the potential $V$ verifies the first condition in Theorem \ref{T: hardytimeindepent} or when we can obtain the additional regularity of the gradient of $u$ in the strip, as in Theorem \ref{T: elcasoantiguo}.
	Here, we choose to prove Theorem \ref{T: hardytimeindepent} using the  Carleman inequality in Lemma \ref{L: lonecesarioparamejorar} in place of the above convexity argument. The reason for our choice is that it is simpler to justify the correctness of the application of the Carleman inequality to a $C([0,1], L^2(\Rn))$ solution to \eqref{E: unaecuacion} than the corresponding monotonicity or  logarithmic convexity  of the solution.

\begin{lemma}\label{L: lonecesarioparamejorar}
The inequality
\begin{multline*}
R\sqrt{\frac{\e}{8\mu}}\,\|e^{\mu |x+Rt(1-t)e_1|^2-\frac{(1+\epsilon)R^2t(1-t)}{16\mu}}g\|_{L^2(\Rm)}\leq\\
\| e^{\mu | x+Rt(1-t)e_1|^2-\frac{(1+\epsilon)R^2t(1-t)}{16\mu}}(\partial_t-i\triangle)g\|_{L^2(\Rm)}
\end{multline*}
holds, when $\e>0$, $\mu >0$, $R>0$ and $g\in C_0^\infty(\Rm)$.
\end{lemma}

\begin{proof}

Let $f=e^{\mu |x+Rt(1-t)|^2-\frac{(1+\epsilon)R^2t(1-t)}{16\mu}}g$. Then,
\begin{equation*}
e^{\mu |x+Rt(1-t)|^2-\frac{(1+\epsilon)R^2t(1-t)}{16\mu}}\left(\partial_t-i\triangle\right)g
=\partial_tf-\mathcal S f- \mathcal A f,
\end{equation*}
and from \eqref{E: loquefcumple}, \eqref{E: formulaoperadores} and \eqref{E: formulaconmutadorindependientetiempo}  with $\gamma =1$, $A+iB=i$ and \[\varphi(x,t)=\mu |x+Rt(1-t)|^2-\tfrac{(1+\epsilon)R^2t(1-t)}{16\mu},\]
 we have
 \begin{equation*}
\begin{split}
&\mathcal S=-4\mu i \left(x+Rt(1-t)e_1\right)\cdot\nabla-2\mu n i\\
&\quad\quad\quad\quad\quad\quad\quad\quad \  +2\mu R(1-2t)\left(x_1+Rt(1-t)\right)-\tfrac{(1+\e)R^2(1-2t)}{16\mu},\\
&\mathcal A=i\triangle+4\mu^2i|x+Rt(1-t)e_1|^2,\\
&\mathcal S_t+[\mathcal S,\mathcal A]=-8\mu\triangle+32\mu^3 |x+Rt(1-t)e_1|^2-4\mu R\left(x_1+Rt(1-t)\right)\\
&\quad\quad\quad\quad\quad\ \ \ +2\mu R^2(1-2t)^2+\tfrac{(1+\e)R^2}{8\mu}-4i\mu R(1-2t)\partial_{x_1}
\end{split}
\end{equation*}
and
\begin{multline}\label{E: unacuenta}
\left(\mathcal S_tf+[\mathcal S,\mathcal A]f,f\right)=32\mu^3 \int |x+Rt(1-t) e_1-\tfrac R{16\mu^2}e_1|^2 |f|^2\,dx+\tfrac{\e R^2}{8\mu}\int |f|^2\,dx\\+8\mu\int |\nabla_{x'} f|^2\, dx
+8\mu\int|i\partial_{x_1}f-\tfrac{R(1-2t)}2f|^2\,dx\ge \tfrac{\e R^2}{8\mu}\int |f|^2\,dx.
\end{multline}

Following the standard method to handle $L^2$-Carleman inequalities \cite{Hor1}, the symmetric and skew-symmetric parts of $\partial_t-\mathcal S-\mathcal A$, as a space-time operator, are respectively $-\mathcal S$ and $\partial_t-\mathcal A$, and its space-time commutator, $[-\mathcal S,\partial_t-\mathcal A]$ is $\mathcal S_t+[\mathcal S,\mathcal A]$. Thus,
\begin{multline}\label{E: loclasico}
\|\partial_tf-\mathcal Sf-\mathcal Af\|^2_{L^2(\Rm)}=\|\partial_tf-\mathcal Af\|^2_{L^2(\Rm)}+ \|\mathcal Sf\|^2_{L^2(\Rm)}\\-2\text{\it Re} \iint \mathcal Sf\overline{\partial_t f-\mathcal Af}\,dxdt\ge\iint[-\mathcal S,\partial_t-\mathcal A]f\overline f\,dxdt
=\int \left(\mathcal S_tf+[\mathcal S, \mathcal A]f, f\right)\,dt,
\end{multline}
and the Lemma \ref{L: lonecesarioparamejorar} follows from \eqref{E: loclasico} and \eqref{E: unacuenta}.
\end{proof}

\begin{proof}[Proof of Theorem \ref{T: hardytimeindepent}]
Let $u$ be as in Theorem \ref{T: hardytimeindepent} and $\widetilde u$, $\widetilde V$ the  corresponding functions defined in Lemma \ref{L: transformada}, when $A+iB=i$. Then, $\widetilde u$ is in $C([0,1], L^2(\Rn))$,
\begin{equation*}
\partial_t\widetilde u=i\left(\triangle \widetilde u +\widetilde V(x,t)\widetilde u\right)\ ,\ \text{in}\  \R^n\times [0,1],
\end{equation*}
$\|e^{\gamma |x|^2}\widetilde u(0)\|$, $\|e^{\gamma |x|^2}\widetilde u(1)\|$ are finite for $\gamma =\frac 1{\alpha\beta}$ and $\gamma >\frac 1{2}$. The proofs of Theorem \ref{T: densityfortimedependentpotentialsvega} or \ref{T: elcasoantiguo} show that in either case
\begin{equation}\label{E: resultadologaritmicosch4}
N_\gamma=\sup_{[0,1]}\|e^{\gamma|x|^2}\widetilde u(t)\|+\|\sqrt{t(1-t)}e^{\gamma|x|^2}\nabla \widetilde u\|_{L^2(\Rn\times [0,1])}<+\infty.
\end{equation}

For given $R>0$, choose $\mu$ and $\e$ such that
\begin{equation}\label{E: escogerparametros}
\frac{(1+\e)^{\frac 32}}{2(1-\e)^3}< \mu\le \frac\gamma{(1+\e)}
\end{equation}
 and let $\theta_M$ and $\eta_R$ be smooth functions verifying, $\theta_M(x)=1$, when $|x|\leq M$, $\theta_M(x)=0$, when $|x|>2M$, $M\ge R$, $\eta_R\in C_0^\infty(0,1)$, $0\le\eta_R\le 1$, $\eta_R(t)=1$ in $[\frac 1R, 1-\frac 1R]$ and $\eta_R=0$ in $[0,\frac 1{2R}]\cup [1-  \frac 1{2R}, 1]$. Then, 
\[g(x,t)=\theta_M(x)\eta_R(t)\widetilde u(x,t),\]
 is compactly supported in $\Rn\times (0,1)$ and
\begin{equation}\label{E: laecuaciondeg}
\partial_tg -i\left(\triangle g+\widetilde Vg\right)=\theta_M\eta_R'\widetilde u-i(2\nabla\theta_M\cdot\nabla \widetilde u+\widetilde u\,\triangle\theta_M )\eta_R.
\end{equation}
The first term on the right hand side of \eqref{E: laecuaciondeg} is supported, where
\begin{equation*}
\mu | x+Rt(1-t)|^2\le \mu(1+\e) |x|^2+\mu\left(1+\frac 1\e\right)\le \gamma |x|^2+ \frac\gamma\e
\end{equation*}
and the second, inside $B_{2M}\setminus B_M\times [\frac 1{2R},1-\frac 1{2R}]$, where 
\[\mu|x+Rt(1-t)e_1|^2\le\gamma |x|^2+\frac{\gamma R^2}\e .\]

Apply now Lemma \ref{L: lonecesarioparamejorar} to $g$ with the values of $\mu$ and $\e$  chosen in \eqref{E: escogerparametros}. This, the bounds for $\mu|x+Rt(1-t)e_1|^2$ in each of the parts of the support of $\partial_tg -i\left(\triangle g+\widetilde Vg\right)$ and the natural bounds for $\nabla \theta_M$, $\triangle\theta_M$ and $\eta_R'$ show that there is a constant $N_\e$ such that  
\begin{equation}\label{E: desigualfinal1}
\begin{split}
&R\|e^{\mu |x+Rt(1-t)e_1|^2-\frac{(1+\epsilon)R^2t(1-t)}{16\mu}}g\|_{L^\infty(\R^n\times [0,1])}\le\\ &N_\e\|\widetilde V\|_{L^\infty(\R^n\times [0,1])}\|e^{\mu |x+Rt(1-t)e_1|^2-\frac{(1+\epsilon)R^2t(1-t)}{16\mu}}g\|_{L^2(\Rn\times [0,1])}\\&+N_\e Re^{\frac\gamma\e}\sup_{[0,1]}\|e^{\gamma |x|^2}\widetilde u(t)\|+N_\e M^{-1} e^{\frac{\gamma R^2}\e}\|e^{\gamma |x|^2}\left(|\widetilde u|+|\nabla\widetilde u|\right)\|_{L^2(\Rn\times [\frac 1{2R},1-\frac 1{2R}])}.
\end{split}
\end{equation}
The first term on the right hand side of \eqref{E: desigualfinal1} can be hidden in the left hand side, when $R\ge 2N_\e\|\widetilde V\|_{L^\infty(\Rn\times [0,1])}$, while the last tends to zero, when  $M$ tends to infinity by \eqref{E: resultadologaritmicosch4}. This and the fact that
$g=\widetilde u$ in $B_{\frac{\e(1-\e)^2R}4}\times [\frac {1-\e}2,\frac {1+\e}2]$, where 
\[\mu |x+Rt(1-t)e_1|^2-\tfrac{(1+\epsilon)R^2t(1-t)}{16\mu}\ge \tfrac{R^2}{16\mu}\left(4\mu^2(1-\e)^6-(1+\e)^3\right)\]
and \eqref{E: escogerparametros} show that
\begin{equation}\label{E: desigualfinal2}
 Re^{C(\gamma,\e)R^2}\|\widetilde u\|_{L^2(B_{\frac R8}\times[\frac{1-\e}2,\frac{1+\e}2] )}\le N_{\gamma, \e}R,
\end{equation}
when $R\ge 2N_\e\|\widetilde V\|_{L^\infty(\Rn\times [0,1])}$. At the same time,
\begin{equation}\label{E: algoinmediatosobreu2}
N^{-1}\|\widetilde u(0)\|\le\|\widetilde u(t)\|\le N\|\widetilde u(0)\|\ ,\ \text{when}\ 0\le t\le 1\ ,\ N= e^{\sup_{[0,1]}\|\Im \widetilde {V}(t)\|_{\infty}},
\end{equation}
and from \eqref{E: resultadologaritmicosch4}
\begin{equation}\label{E: concentracion}
\|\widetilde u(t)\|\le \|\widetilde u(t)\|_{L^2(B_{\frac R8})}+e^{-\frac{\gamma R^2}{64}}N_\gamma\ ,\ \text{when} \ 0\le t\le 1.
\end{equation}
Then,  \eqref{E: desigualfinal2}, \eqref{E: algoinmediatosobreu2} and \eqref{E: concentracion} show that there is a constant $N_{\gamma, \e, V}\,$, which depends on $N_\gamma$, $\e$ and the $\sup_{[0,1]}\|V(t)\|_\infty$, such that
 \begin{equation*}
e^{C(\gamma,\e)R^2}\|\widetilde u(0)\|\le N_{\gamma, \e, V}.
\end{equation*}
Let then $R$ tend to infinity to derive that $u\equiv 0$.
\end{proof}
\end{section}
\begin{section}{A positive commutator and a misleading frequency function}\label{S:  A positive commutator and a misleading frequency function}
When $f=e^{a(t)|x|^2}u$ and $u$ is solution to the free Schr\"odinger equation  in $\R\times [-1,1]$, $f$ verifies $\partial_tf=\mathcal Sf+\mathcal Af$,
with symmetric and skew-symmetric operators,
\begin{equation*}
\mathcal S=-4ia\left(x\partial_x+\tfrac 12\right)+a'x^2\quad,\quad\mathcal A= i\left(\partial^2_x+4a^2x^2\right).
\end{equation*}
In this case (See \eqref{E: formulaconmutadorindependientetiempo})
\begin{equation*}
\mathcal S_t+\left[\mathcal S,\mathcal A\right]= \tfrac{2a'}a\mathcal S-8a\partial^2_x+\left(32a^3+a''-\tfrac{2a'^2}a\right)x^2,
\end{equation*}
and if $a$ is a positive and even solution of
\begin{equation}\label{E: la ecuaci—nrara}
32a^3+a''-\tfrac{2a'^2}a=0,\ \text{in}\ [-1,1],
\end{equation}
the formal calculations in Lemma \ref{L:  freq1}  show that $H_a(t)=\|e^{a(t)x^2}u(t)\|^2$ verifies
\begin{equation*}
\partial_t\left(a^{-2}\partial_t\log{H_a}(t)\right)\ge 0,\ \text{in}\ [-1,1],
\end{equation*}
and the integration of the inequality 
\begin{equation*}
a^2(\tau)\partial_s\log{H_a(s)}\le a^2(s)\partial_\tau\log{H_a(s)}\ ,\ \text{when}\ -1\le s\le 0\le\tau\le 1,
\end{equation*}
implies that
\begin{equation}\label{E:  convexidadlogaritmicaenotrointervalo}
H_a(0)\le H_a(-1)^{\frac 12}H_a(1)^{\frac 12}.
\end{equation}
On the other hand, if $a$ solves
\begin{equation*}
\begin{cases}
32a^3+a''-\frac{2a'^2}{a}=0,\\
a(0)=1,\, a'(0)=0,
\end{cases}
\end{equation*}
$a$ is positive, even and $\lim_{R\rightarrow+\infty} Ra(R)=0$. Moreover, $a_R(t)=Ra(Rt)$ also  solves \eqref{E: la ecuaci—nrara}, and if the formal calculation is correct for $H_{a_R}$, 
\eqref{E:  convexidadlogaritmicaenotrointervalo} would imply that
\begin{equation*}
\|e^{Rx^2}u(0)\|^2\le \|e^{Ra(R)x^2}u(-1)\|\|e^{Ra(R)x^2}u(1)\| .
\end{equation*} 
In particular, $u\equiv 0$; but 
\[u(x,t)=\left(t-i\right)^{-\frac 12}e^{\frac{i|x|^2}{4(t-i)}}\]
contradicts this

This shows that there are functions $\varphi$, which make non-negative the commutator of the symmetric and skew-symmetric parts of $e^{\varphi}\left(\partial_t-i\partial^2_x\right)e^{-\varphi}$ and such that it is not possible to plug in or enter  in the associated Carleman inequality or frequency function some reasonable solutions of the free Schr\"odinger equation. It also shows that the rather complex arguments  we used to derive the  logarithmic convexity of  
 \begin{equation*}
H(t)=\int_{\Rn}e^{2\gamma |x|^2}|u(t)|^2\,dx ,
\end{equation*}
are in fact necessary, when $u$ in $C([0,1],L^2(\Rn))$ is a solution verifying the conditions in Lemma \ref{L: convexlogsch} and as long as a more suitable representation formula for these solutions is not available. By suitable we mean a formula  which allows to derive the quadratic exponential decay of the solution in the interior of a time slab from the known decay of the solution at the top and bottom of the slab.
\end{section}
\begin{section}{Parabolic analog. Proof of Theorem \ref{T: toremaparabolicco}}\label{S: Parabolic analog}

Assume that $u$ verifies the conditions in Theorem \ref{T: toremaparabolicco} and let $\widetilde u$ be the conformal or Appel transformation of $u$ defined in Lemma \ref{L: transformada} with $A+iB=1$, $\alpha=1$ and $\beta=1+\frac 2\delta$. $\widetilde u$ is in $L^\infty([0,1]), L^2(\Rn))\cap L^2([0,T], H^1(\Rn))$, verifies
\[
\partial_t\widetilde u=\triangle\widetilde u+\widetilde V(x,t)\widetilde u,\ \text{in}\ \Rn\times (0,1]
\]
with $\widetilde V$ a bounded potential in $\Rn\times [0,1]$ and if $\gamma=\frac 1{2\delta}$, we have
\begin{equation*}
\|e^{\gamma |x|^2}\widetilde u(0)\| =\|u(0)\|\ ,\ \|e^{\gamma |x|^2}\widetilde u(1)\|=\|e^{\frac{|x|^2}{\delta^2}}u(1)\|
\end{equation*}

 From Lemma \ref{L: convexlogsch} and Lemma \ref{L: regularidad} with $A+iB=1$, we have
\begin{multline}\label{E: regularidad3}
\sup_{[0,1]}\|e^{\gamma |x|^2}\widetilde u(t)\|+\|\sqrt{t(1-t)}e^{\gamma |x|^2}\nabla\widetilde u\|_{L^2(\Rn\times [0,1])}\le\\ e^{N(M_1+M_1^2)}\left(\|e^{\gamma |x|^2}\widetilde u(0)\|+ \|e^{\gamma |x|^2}\widetilde u(1)\|\right)
\end{multline}
where $M_1=\|\widetilde V\|_{L^\infty(\Rn\times [0,1])}$. The proof is finished by pluging in 
\[g(x,t)=\theta_M(x)\eta_R(t)\widetilde u(x,t)\]
in the Carleman inequality below and in complete analogy with the argument we used   to prove Theorem \ref{T: hardytimeindepent}.

\begin{lemma}\label{L: carlemanparabolico}
The inequality
\begin{multline*}
R\sqrt{\frac{\e}{8\mu}}\,\|e^{\mu |x+Rt(1-t)e_1|^2+\frac{R^2t(1-t)(1-2t)}6-\frac{(1+\epsilon)R^2t(1-t)}{16\mu}}g\|_{L^2(\Rm)}\leq\\
\| e^{\mu |x+Rt(1-t)e_1|^2+\frac{R^2t(1-t)(1-2t)}6-\frac{(1+\epsilon)R^2t(1-t)}{16\mu}}(\partial_t-\triangle)g\|_{L^2(\Rm)}
\end{multline*}
holds, when $\e>0$, $\mu >0$, $R>0$ and $g\in C_0^\infty(\Rm)$.
\end{lemma}

\begin{proof}
Let $f=e^{\mu |x+Rt(1-t)e_1|^2+\frac{R^2t(1-t)(1-2t)}6-\frac{(1+\epsilon)R^2t(1-t)}{16\mu}}g$. Then,
\begin{equation*}
e^{\mu |x+Rt(1-t)e_1|^2+\frac{R^2t(1-t)(1-2t)}6-\tfrac{(1+\epsilon)R^2t(1-t)}{16\mu}}\left(\partial_t-\triangle\right)g
=\partial_tf-\mathcal S f- \mathcal A f,
\end{equation*}
and from \eqref{E: loquefcumple}, \eqref{E: formulaoperadores} and \eqref{E: formulaconmutadorindependientetiempo}  with $\gamma =1$, $A+iB=1$ and 
\begin{equation*}
\varphi(x,t)=\mu |x+Rt(1-t)e_1|^2+\tfrac{R^2t(1-t)(1-2t)}6-\tfrac{(1+\epsilon)R^2t(1-t)}{16\mu}\ ,
\end{equation*}
 we have
 \begin{equation*}
\begin{split}
&\mathcal S=\triangle+4\mu^2|x+Rt(1-t)e_1|^2+2\mu R(1-2t)\left(x_1+Rt(1-t)\right)\\
&\quad+(t^2-t+\tfrac 16)R^2-\tfrac{(1+\e)R^2(1-2t)}{16\mu},\\
&\mathcal A=-4\mu \left(x+Rt(1-t)e_1\right)\cdot\nabla-2\mu n,\\
&\mathcal S_t+[\mathcal S,\mathcal A]=-8\mu\triangle+32\mu^3 |x+Rt(1-t)e_1|^2+2\mu R^2(1-2t)^2\\
&+4\mu R(4\mu (1-2t)-1)(x_1+Rt(1-t))+(2t-1)R^2+\tfrac{(1+\e)R^2}{8\mu}
\end{split}
\end{equation*}
and
\begin{multline}\label{E: unacuenta1}
\left(\mathcal S_tf+[\mathcal S,\mathcal A]f,f\right)=32\mu^3 \int |x+Rt(1-t) e_1+\tfrac {(4\mu(1-2t)-1)R}{16\mu^2}e_1|^2 |f|^2\,dx\\+8\mu\int |\nabla f|^2\, dx++\tfrac{\e R^2}{8\mu}\int |f|^2\,dx\ge \tfrac{\e R^2}{8\mu}\int |f|^2\,dx.
\end{multline}
Finally,
\begin{equation*}
\varphi(x,\tfrac 12)=\mu |x+\tfrac R4e_1|^2-\tfrac{(1+\epsilon)R^2}{64\mu}\ge (4\mu^2(1-\e)^2-(1+\e))\tfrac{R^2}{64\mu},
\end{equation*}
when $|x|\le \frac{\e R}4$, and it is positive for $\mu>\tfrac 12$ and $\e>0$ small.
\end{proof}
\begin{remark}\label{R: 2}
\eqref{E: regularidad3}, \eqref{E: unacuenta1} and the interior regularity of parabolic equations show that the formal calculations in Lemma \ref{L: freq1} to prove the logarithmic convexity of
\begin{equation*}
H_\mu(t)=\int_{\Rn}e^{2\mu |x+Rt(1-t)e_1|^2+\frac{R^2t(1-t)(1-2t)}3-\frac{R^2t(1-t)}{8\mu}}|\widetilde u(t)|\,dx
\end{equation*}
are correct, when $\mu<\gamma$. In particular,
\begin{equation*}
H_\mu(\tfrac 12)\le e^{N(M_1+M_1^2)}H_\mu(0)^{\frac 12}H_\mu(1)^{\frac 12},
\end{equation*}
and letting $\mu$ increase to $\gamma$ and then $R$ tend to infinity, one also gets that $\widetilde u\equiv 0$ in $\Rn\times [0,1]$, when $\gamma>\tfrac 12$.
\end{remark}
\end{section}



\begin{thebibliography}{99}



\bibitem{bonamie}  A. Bonami, B. Demange, A survey on uncertainty principles related to quadratic forms. Collect. Math. 2006, Vol. Extra, 1--36. 

\bibitem{ekpv06} L. Escauriaza, C.E. Kenig, G. Ponce, L. Vega,  \emph{On Uniqueness  Properties of Solutions of Schr\"odinger  Equations,} Comm. PDE. {\bf 31}, 12 (2006) 1811--1823.

\bibitem{ekpv07} L. Escauriaza, C.E. Kenig, G. Ponce, L. Vega,  \emph{On  Uniqueness Properties of Solutions of the k-generalized KdV,} J. of Funct. Anal. {\bf 244}, 2 (2007) 504--535.


\bibitem{ekpv08} L. Escauriaza, C.E. Kenig, G. Ponce, L. Vega,  \emph{Convexity of Free Solutions of Schr\"odinger Equations with Gaussian Decay.} To appear.

\bibitem{e98} L.C. Evans,  \emph{Partial Differential Equations.} Amer. Math. Soc. (1998)

\bibitem{Hor1} L.  H\"ormander, \emph{Linear partial differential operators}, Berlin, Springer (1969).

\bibitem{Ioke04} A. D. Ionescu, C. E. Kenig, \emph{$L^p$-Carleman inequalities and uniqueness of solutions of nonlinear Schr\"odinger equations,} Acta Math. {\bf 193}, 2 (2004) 193--239.

\bibitem{Ioke06} A. D. Ionescu, C. E. Kenig, \emph{Uniqueness properties of solutions of Schr\"odinger equations,} J. Funct. Anal. {\bf 232} (2006) 90--136.

\bibitem{Isakov93} V. Isakov, \emph{Carleman type estimates in anisotropic case and applications, } J. Diff. Eqs. {\bf 105} (1993) 217--238.

\bibitem{j06} P. Jaming, \emph{Uncertainty  Principles for Orthonormal Bases.} arxiv.org/pdf/math/0606396

\bibitem{kpv02} C.E. Kenig, G. Ponce, L. Vega, \emph{On unique continuation for nonlinear Schr\"odinger equations}, Comm. Pure Appl. Math. {\bf 60} (2002) 1247--1262.

\bibitem{kpv03}  C.E. Kenig, G. Ponce, L. Vega, \emph{On unique continuation of solutions to the generalized KdV equation},  Math. Res. Letters {\bf10} (2003) 833--846.

\bibitem{lm60}  J.L. Lions, B. Malgrange, \emph{Sur l'unicit e r etrograde dans les probl emes mixtes paraboliques,} Math. Scan. {\bf8} (1960) 277--286.

\bibitem{pazy}  A. Pazy, \emph{Semigroups of linear operators with application to partial differential equations}.  Springer-Verlag, Berlin, New York (1983).


\bibitem{SiSu} A. Sitaram, M. Sundari, S. Thangavelu, {Uncertainty principles on certain Lie groups},  Proc. Indian Acad. Sci. Math. Sci. {\bf 105} (1995), 135-151

\bibitem{StSh} E.M. Stein, R. Shakarchi, \emph{Princeton Lecture in Analysis II. Complex Analysis,} Princeton University Press.


\bibitem{Than} S. Thangavelu, \emph{Lectures on Hermite and Laguerre expansions}, Princeton Univ. Press, Princeton, New Jersey (1993).



\end{thebibliography}
\end{document}